\documentclass[preprint,10pt]{elsarticle}

\usepackage{amssymb}
\usepackage{amsfonts, amsmath}
\usepackage{graphicx}
\usepackage{tikz}
\graphicspath{ {./images/}}
\usepackage{caption}
\usepackage{subcaption}
\usepackage{amssymb}
\usepackage{amsthm}
\usepackage{float}

\newtheorem{theorem}{Theorem}[section]
\newtheorem{lemma}[theorem]{Lemma}
\newtheorem*{remark}{Remark}
\newtheorem{proposition}[theorem]{Proposition}

\usepackage{varwidth}
\usepackage[linesnumbered,ruled,vlined]{algorithm2e}
\DeclareCaptionFormat{myformat}{
  \begin{varwidth}{\linewidth}%
    \centering
    #1#2#3%
  \end{varwidth}%
}

\usepackage{geometry} 
\journal{Journal of Computational and Applied Mathematics}

\begin{document}

\begin{frontmatter}



\title{Uniform Approximation of Functions with Asymmetric Growth and Decay by Deep Weighted Polynomials}


\author[uchic]{Kingsley Yeon}

\affiliation[uchic]{organization={Department of Statistics and CCAM, University of Chicago},
            city={Chicago},
            postcode={60637}, 
            state={IL},
            country={USA}}

\author[zbmath]{Steven B. Damelin}

\affiliation[zbmath]{organization={Department of Mathematics, ZBMATH-OPEN},
            addressline={Hermann-von-Helmholtz-Platz 1}, 
            city={Eggenstein-Leopoldshafen},
            postcode={76344}, 
            country={Germany}}

\begin{abstract}
Functions that grow without bound on one side of the real line and decay to zero on the other cannot be approximated uniformly by ordinary polynomials on unbounded domains. Motivated by classical weighted polynomial approximation, we introduce a class of one-sided weighted \emph{deep} (composite) polynomial approximants for such asymmetric targets. The weight suppresses polynomial growth on the decaying side, while the composite polynomial remains free to capture growth on the other side. We prove that this mechanism reduces the half-line approximation problem to approximation on a compact interval whose length grows slowly with the degree, and we establish density and existence of best approximants in the appropriate closure of the model class. For computation, we first formulate the method as a trainable computational graph for \emph{deep} weighted polynomial approximation. However, direct end-to-end optimization becomes increasingly ill-conditioned at high composite degree and can suffer from local minima. To address this, we introduce a fine-tuning procedure in which a fixed inner composition of monotone polynomial self-maps supplies the effective degree, while only the outer polynomial and weight parameters are trained; the outer fit reduces to a linear program. Numerical experiments on Black--Scholes option-pricing functions show that the resulting fine-tuned weighted \emph{deep} polynomial achieves smaller uniform and \(L_2\) errors than matched-budget polynomial baselines and resolves the decaying tail to machine precision.
\end{abstract}




\end{frontmatter}



\section{Introduction}\label{sec:intro}

Let \(\{p_\ell\}_{\ell}^L\) be a family of univariate polynomials on the real line where each \(p_\ell\) has degree \(d_\ell\). Here $d_{\ell}\geq 1$. A \emph{deep polynomial} of \(L\) layers is defined via the composition
\[
P(x):=p_L(p_{L-1}(p_{L-2}(...(p_1(x)))).
\]
The total (composite) degree of \(P\) is
\[
D \;=\; d_1\,d_2\,\cdots\,d_L.
\]
Typically, the number of layers $L$ is small in comparison to the total degree $D$ which can be taken to be quite large. 
The composition $P$ is parametrized by the coefficients of its $L$ layers. Without any normalization, these coefficients are not unique: at each of the $L-1$ interfaces between consecutive layers, an invertible affine reparametrization $p_\ell\mapsto \beta\,p_\ell+\gamma$ ($\beta\neq 0$) can be absorbed into the adjacent layer without changing $P$, so the representation is redundant and the optimization landscape contains flat directions. Following \cite{KY}, we therefore \emph{normalize} each inner layer ($\ell < L$) to be monic with zero constant term, which removes this affine ambiguity, makes the representation unique, and renders the gradient-based optimization substantially more stable. The raw coefficient count of $P$ before normalization is $\sum_{\ell=1}^L(d_\ell+1)=\sum_{\ell=1}^L d_\ell+L$; each of the $L-1$ inner interfaces carries a two-parameter affine redundancy (one scale and one shift), and eliminating monic leading coefficient and vanishing constant term in each of the $L-1$ inner layers removes exactly $2(L-1)$ of these coefficients. Hence after normalization the number of genuinely free parameters is
\[
n_{\mathrm{deep}} \;=\; \Bigl(\sum_{\ell=1}^{L}d_\ell\Bigr) + L - 2(L-1)\;=\;\Bigl(\sum_{\ell=1}^{L}d_\ell\Bigr) - L + 2,
\]
and it is this count that we use when matching the parameter budget between the deep polynomial and a classical baseline of the same number of free coefficients. It is widely believed that the layered structure of composite polynomials mimics in many ways the layers of neurons in deep learning networks, and hence their approximation properties are natural to study \cite{Dev}.

Functions with singularities are often difficult to approximate by polynomials due to poor behavior near the singularity, which is why classical theory relies on smoothness. Composite polynomials, however, often capture such structure more effectively. At its core, approximation theory is a story of compression: we seek low-dimensional parameterizations that accurately represent a function. Polynomials excel at compressing smooth functions, with classical orthogonal polynomials achieving exponential convergence for analytic functions. Similarly, Fourier series provide optimal representations for periodic inputs in $L^2$. For continuous functions with sharp transitions or steep gradients, rational functions are often far more efficient than polynomials; their best (minimax) approximants can be computed by the Remez algorithm, which applies to continuous targets, and near-best rational approximants can be computed in barycentric form by the AAA algorithm \cite{AAA}.

In practice, many numerical algorithms use structured iterations such as Newton's method. These can be written as composite rational or polynomial updates \( f_{k+1}(x) = g(f_k(x), x) \), where \( g \in \mathcal{P}_n \) is a fixed-degree polynomial. Then \( f_k \) becomes a deep polynomial of degree \( n^k \), defined by few tunable parameters. Such iterative structures appear naturally in algorithms for non-arithmetic operations. For example, the sign function can be written as \( \text{sgn}(x) = x / |x| \), the comparison operator as \( \text{comp}(u, v) = \tfrac{1}{2}(\text{sgn}(u - v) + 1) \), and the max operation as \( \max(u, v) = \tfrac{1}{2}(u + v + (u - v)\text{sgn}(u - v)) \). Newton's method is often employed to approximate \( \text{sgn}(x) \), and has been shown to be optimal in the sense of requiring the fewest non-scalar multiplications \cite{Lee}.

As a concrete example, to approximate \( |x| \), consider
\[
f_{k+1}(x) = \tfrac{1}{2} f_k(x) \left(3 - x^2 f_k(x)^2\right), \quad f_0(x) = 1,
\]
which converges to \( 1/|x| \) for \( x \ne 0 \). Multiplying by \( x^2 \) yields a sequence of deep polynomials converging exponentially to \( |x| \) \cite{KY}. In many cases, such composite iterations surpass classical minimax approximants while using fewer degrees of freedom.

Our work is motivated by \cite{KY} and \cite{GN}. Classically, $|x|$ admits type-$(d,d)$ rational minimax approximation on $[-1,1]$ with rate $\exp(-\sqrt{d})$ \cite{Sta}, and for fixed $p\geq 2$ the function $x^{1/p}$ admits an explicit $\lceil\log_{p}d\rceil+1$ layer composite rational approximant on $[0,1]$ converging double exponentially in its number of degrees of freedom \cite{GN}, whereas Jackson's theorem limits polynomials of degree $d$ to the algebraic rate $1/d$ for these targets \cite{Lub2}. Deep polynomials recover exponential rates in the number of free parameters for the same targets \cite{KY}.

Motivated by the property that deep polynomials are well-suited for approximating non-regular functions (e.g., those that are only \(C^0\)), we ask whether they can be combined with a one-sided decaying weight to approximate functions that grow unbounded on one side and decay on the other. A simple example is the exponential function. Polynomials alone cannot approximate \(e^{-x}\) uniformly on \(\mathbb{R}\), since they diverge at infinity. Rational functions do converge, but require divisions. Instead, we consider weighted polynomials of the form \(P(x)\,e^{-Q(x)}\), where the external field satisfies \(Q(x)=0\) for \(x<0\) (so the weight is \(1\) there) and \(Q(x)\ge 0\) for \(x\ge 0\). For example, taking \(Q(x)=x^2\) on \(x\ge0\) gives a one-sided Gaussian weight. On the decaying side \(x\ge0\) the weight suppresses the polynomial's growth so that, by the restricted-range mechanism described in Section~\ref{sec:restricted}, the weighted approximant is effectively controlled out to large $x$; on the growing side \(x<0\) the weight is inactive, and the composite polynomial $P$ must reproduce the growth $e^{|x|}$ on a finite interval, which is ordinary (unweighted) polynomial approximation of a smooth function. The content of the construction is that a single low-parameter weighted composite resolves both regimes at once: the weight reduces the unbounded decaying side to a compact effective interval, made precise in Proposition~\ref{prop:compact}, while the composite polynomial supplies the resolution required on the growing side. The weighted product is only $C^1$ at the origin, which makes a deep polynomial a natural choice for the resulting non-smooth interface.

\begin{figure}[H]
    \centering
    \includegraphics[width=0.6\textwidth]{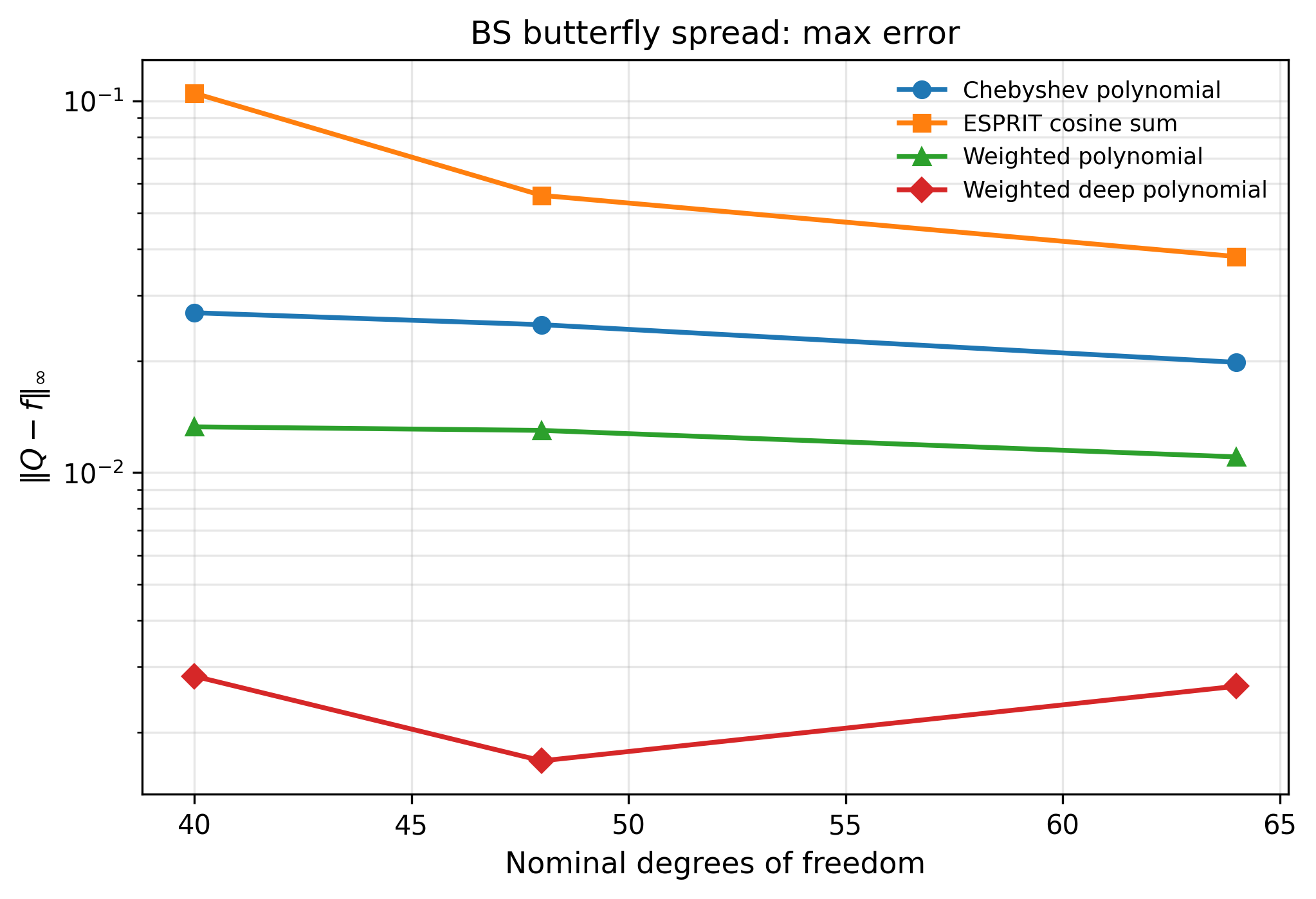}
    \caption{Black--Scholes butterfly benchmark comparing the cosine-sum ESPRIT baseline from Algorithm~1 of \cite{DPR-ESPRIT}, the Chebyshev polynomial, the learned weighted polynomial, and the fine-tuned weighted deep polynomial. At the matched budgets shown, ESPRIT gives the largest uniform error, while the fine-tuned weighted deep model gives the smallest.}
    \label{fig:bs_butterfly_esprit}
\end{figure}

\begin{figure}[H]
    \centering
    \includegraphics[width=0.8\textwidth]{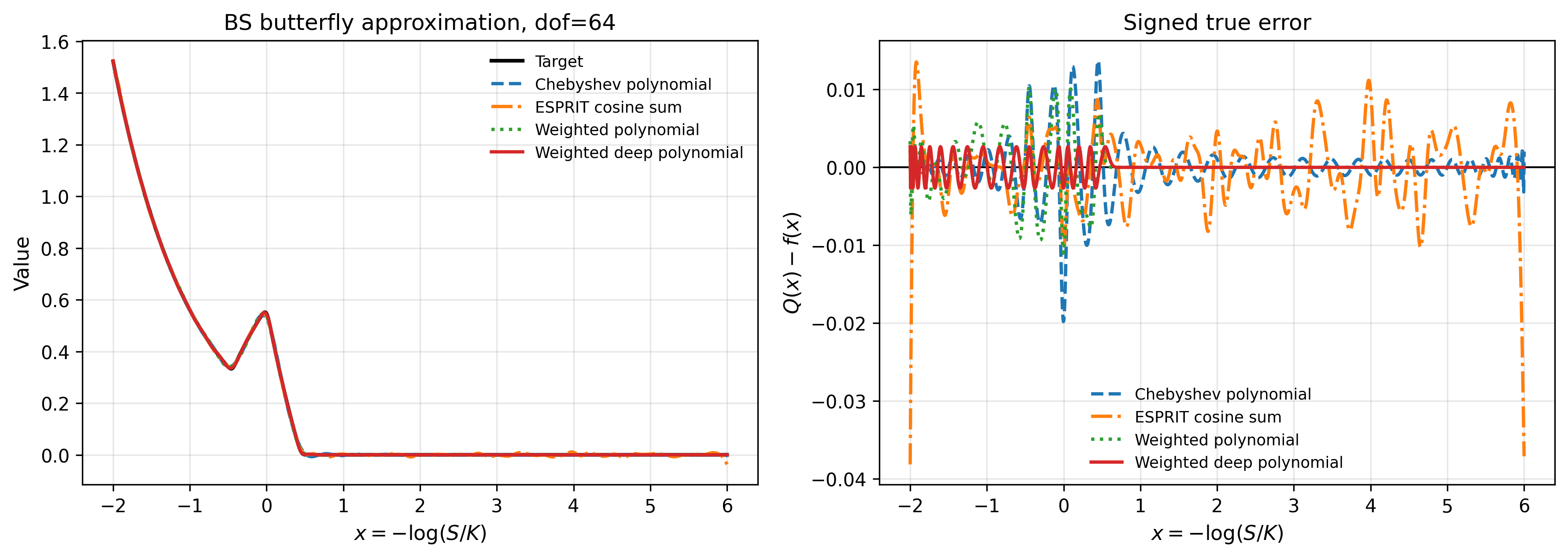}
    \caption{Black--Scholes butterfly approximation and signed error at the largest budget. The ESPRIT cosine sum leaves a visibly larger global residual, while the weighted deep polynomial concentrates the residual near the localized transition.}
    \label{fig:bs_butterfly_esprit_approx}
\end{figure}

\begin{figure}[H]
    \centering
    \includegraphics[width=0.6\textwidth]{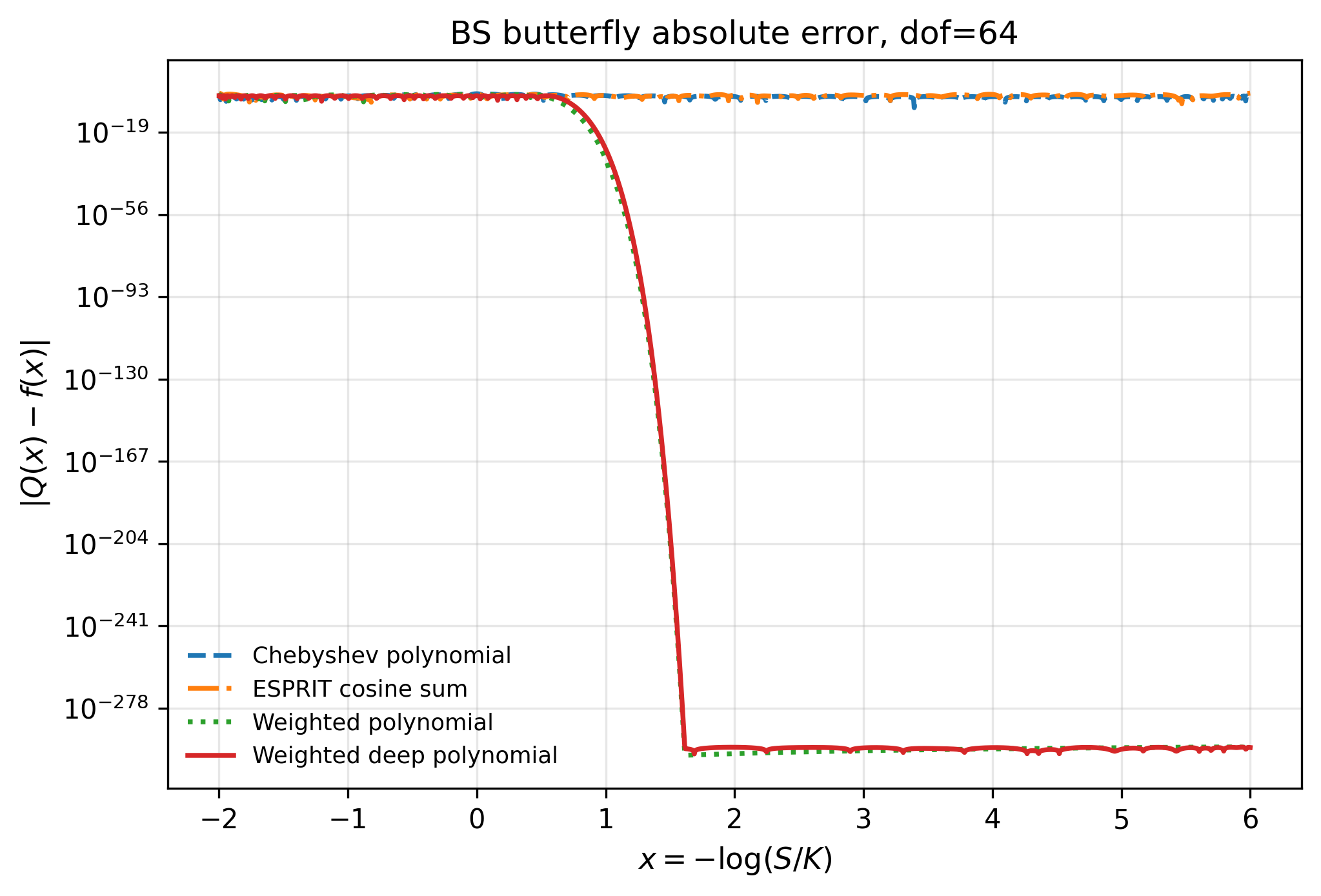}
    \caption{Black--Scholes butterfly pointwise absolute error on a logarithmic scale. The ESPRIT fit remains the least accurate across the interval, whereas the weighted deep polynomial damps the decaying side through the learned weight and achieves substantially smaller error near the relevant transition region.}
    \label{fig:bs_butterfly_esprit_detail}
\end{figure}

A natural nonlinear alternative for such tasks is approximation by exponential or cosine sums. We emphasize that \emph{pure} exponential sums $\sum_n a_n e^{x_n t}$ are not dense in $C[a,b]$; density requires the \emph{extended} or confluent family $\sum_n\sum_i a_{ni} t^i e^{x_n t}$, in which coalescing frequencies generate the polynomial prefactors $t^i$ \cite{Braess}. Fitting these models is nonlinear in the frequencies. Classical Prony-type and subspace methods, including ESPRIT, estimate the frequencies through a matrix-pencil or generalized eigenvalue problem and then recover the amplitudes by linear least squares \cite{ESPRIT,DPR-ESPRIT}. This makes ESPRIT an attractive comparison point for short smooth oscillatory signals, but our Black--Scholes butterfly experiment in Figures~\ref{fig:bs_butterfly_esprit}--\ref{fig:bs_butterfly_esprit_detail} shows that the cosine-sum ESPRIT baseline is the least accurate method at the matched budgets considered: its uniform error remains at the $10^{-1}$ scale, whereas the fine-tuned weighted deep polynomial reaches the $10^{-3}$ scale. This is not surprising, since the experiment is an approximation problem rather than an exact sparse-recovery problem. ESPRIT is designed to recover a short real cosine sum from midpoint samples; for a localized, non-periodic, one-sided decaying option profile, the Toeplitz--Hankel pencil must infer global frequencies from data that are nearly flat on much of the interval but sharply localized near the strikes. Deep weighted polynomials are a division-free alternative. The one-sided weight models the decay, and the stable inner map reallocates degrees of freedom toward the region where the payoff changes sharply. This structure is especially relevant for functions such as \(\log x\), which exhibit asymmetric behavior and arise frequently in applications such as computing \(\log\det\) \cite{KY1}.

\subsection{Contributions and outline}

The contributions of this paper are the following. First, we formulate one-sided weighted deep polynomial approximation for targets with asymmetric growth and decay, and we prove a quantitative reduction of the half-line problem to a compact interval (Proposition~\ref{prop:compact}): the weight makes the problem effectively compact, with an effective interval of length $O\bigl((n\log n)^{1/\beta}\bigr)$ for a degree-$n$ approximant and field exponent $\beta$. Sections~\ref{sec:bernstein} and \ref{sec:jackson} survey the classical restricted-range and weighted Jackson--Bernstein theory on which this reduction rests; the expert reader may proceed directly to Section~\ref{sec:weighteddeep}. Second, we prove density of weighted deep polynomials in $C(K)$ and existence of best approximants in the closure of each fixed architecture, identifying the lack of closedness responsible for possible non-attainment within the architecture itself (Section~\ref{sec:existence}). Third, we give a computational-graph algorithm for training all layers (Section~\ref{sec:computational}) and, to address its ill-conditioning at large depth, a fine-tuning procedure in which the inner composition is fixed and only the outer polynomial and weight are trained by a convex fit (Section~\ref{sec:finetune}). Fourth, we validate the framework on a suite of Black--Scholes option-pricing functions at matched parameter budgets and evaluation cost.

\subsection{Notations}

We adopt the following conventions throughout. $a,b,C,C_1, ...$ will denote positive finite constants. Their dependence on different quantities will be indicated and the same symbol may denote a different constant at any given time. The context will be made clear.
For positive sequences $b_n,c_n$, we write $b_n=O(c_n)$ if uniformly in $n$, $b_n\leq Cc_n$. We also write
$a_n\sim b_n$ if $b_n=O(c_n)$ and $c_n=O(b_n)$ uniformly in $n$.
Similar notations will be used for functions. When we write 
$||.||_{.}$ we mean the $L_{\infty}$ norm. Unless stated otherwise, by a polynomial we mean an algebraic polynomial.

\section{Bernstein's approximation problem}\label{sec:bernstein}
In 1924, the mathematician Sergei Bernstein \cite{B} came up with a problem known as ``Bernstein's approximation problem''. In analogy to the Weierstrass approximation theorem which tells us that polynomials are dense in the space of real valued continuous functions on a compact interval, Bernstein asked if there are analogues for the real line? This question has ramifications which still continue to be studied. The first immediate observation is that we have to deal with the unboundedness of polynomials on the real line. Thus the formulation of Bernstein's problem includes the combination $Sw$, where the function $w$ decays fast at $\pm \infty$ to counter the growth of the polynomial say $S$. We call the function $w$ a weight. We require from $w$ that 
$\lim_{|x|\to \infty}x^nw(x)=0,\, n\geq 0$. What can be approximated and in what sense?
\medskip

A precise statement of Bernstein's problem is the following. When is it true that for every continuous 
$f:\mathbb R\to \mathbb R$ with $\lim_{|x|\to \infty}(fw)(x) =0$, there exists a sequence of 
polynomials $\left\{S_n\right\}_{n=1}^{\infty}$ with 
$\lim_{n\to \infty}||(f-S)w||_{\mathbb R}=0$? If true, then we say that Bernstein's problem has a positive solution. The restriction that the combination $fw$ has limit zero for large argument is essential. Indeed, if $x^nw(x)$ is bounded on the real line for every non-negative integer $n$, then $x^nw(x)$ has limit zero for large argument and so the same is true for every weighted polynomial $Sw$. So we could not hope to approximate, in the uniform norm a function $f$ for which $fw$ does not have limit $0$ for large argument.

Bernstein's problem has been studied, solved and extended by many, including for example N.I Akhiezer, K.I Bambenko, L.de Branges, L. Carleson, M.M Dzbasjan, T. Hall, S.Izumi, T. Kawata, S.N Mergelyan and H. Pollard and V.S Videnskii. See for example \cite{Koo} and the many references cited therein. 
\medskip

One formulation of an answer to Bernstein's problem is the following: 
\begin{itemize}
\item[(1)] Let $w:\mathbb R\to (0,1]$ be continuous with 
$\lim_{|x|\to \infty}x^nw(x)=0,\, n=0,1,2,...$. Then Bernstein's problem has a positive answer iff (a) and (b) below hold:
\begin{itemize}
\item[(a)] $\int_{\mathbb R}\frac{\log(w^{-1}(t))}{1+t^2}dt=\infty.$
\item[(b)] There exists a sequence of polynomials $\{S_n\}$ of degree at most $n\geq 1$ such that $\lim_{|x|\to \infty}(S_nw)(x)=0$ and with ${\rm sup}_{n\geq 1}||S_nw||_{\mathbb R}<\infty$.
\end{itemize}
\end{itemize}
In particular:
\begin{itemize}
\item[(2)] Let $w:\mathbb R\to (0,1]$ be even, with $\log(w^{-1}(e^x))$ convex
on $\mathbb R$
and with $\lim_{|x|\to \infty}x^nw(x)=0,\, n=0,1,2,...$. Then Bernstein's problem has a positive answer iff $\int_{0}^{\infty}\frac{\log(w^{-1}(t))}{1+t^2}dt=\infty.$ 
\end{itemize}

An immediate observation regarding the condition $\int_{0}^{\infty}\frac{\log(w^{-1}(t))}{1+t^2}dt=\infty$ is that
it forces the weight $w$ to decrease at least as fast as a polynomial for large argument. 
\medskip

Here is one example:
\medskip

Let $\lambda\geq 1$ and $x\in \mathbb R$ and let $w_{\lambda}(x):=\exp(-|x|^{\lambda})$. In the literature, these weights (the special case $\lambda=2$ is just the Hermite weight)
are examples of what are called Freud weights after the mathematician 
Geza Freud \cite{Mha} who first investigated them. Here, $w_{\lambda}$ is of smooth polynomial decay for large argument and Freud weights are essentially characterized by the property that they are of smooth polynomial decay for large argument.

\section{Rates of convergence: weighted Jackson--Bernstein theorems}\label{sec:jackson}

In this section, we look at rates of convergence for the Bernstein's density theorem. We will focus on Freud weights but these examples generalize to more general classes of weights. See for example 
\cite{Lub, Dam, Dam1, Dam2, Lub1, Mha} and the references cited therein.
\medskip

\subsection{ The classical case} 

The classical Jackson approximation theorem tells us, for example, that on a finite real interval $[a,b]$, for every absolutely continuous function $f:[a,b]\to\mathbb R$ with $f'\in L^\infty[a,b]$ one has, with a constant $C$ independent of $n$ and $f$,
\[
E_n[f]={\rm inf}_{{\deg}S\leq n}\,\|f-S\|_{[a,b]}\leq C
\left(\frac{b-a}{n}\right)\|f'\|_{[a,b]},\qquad n\geq 1 .
\]
Here the ${\rm inf}$ is taken over the $(n+1)$-dimensional space of algebraic polynomials of degree at most $n$, and is attained by standard compactness arguments. The rate is best possible among absolutely continuous, real-valued functions $f:[a,b]\to\mathbb R$ whose derivative is bounded. More generally, for any continuous $f$ with modulus of continuity
\[
\omega(f,\delta):={\rm sup}\left\{|f(x)-f(y)|:\, x,y\in [a,b]\ {\rm and}\ |x-y|\leq \delta\right\},\quad \delta>0,
\]
one has, again with $C$ independent of $n$ and $f$,
\[
E_n[f]\leq C\,\omega\!\left(f,\tfrac{b-a}{n}\right),\qquad n\geq 1 .
\]
(The right-hand side is finite for every continuous $f$ on the compact interval $[a,b]$, so no extra hypothesis is needed here.)
\medskip

If $f:[a,b]\to\mathbb R$ has a bounded $k$th derivative ($k\geq 1$), the minimax rate above improves to
\[
E_n[f]\leq C_k\,\frac{(b-a)^k}{n^k}\,\|f^{(k)}\|_{[a,b]},\qquad n\geq 1,
\]
where the constant $C_k$ depends only on $k$ and is independent of $n$. (The bound is uniform in $n$; it is of course not uniform in $f$, since it scales with $\|f^{(k)}\|_{[a,b]}$, nor in $k$.)
\medskip

Converse Bernstein results are also classical in the following sense. We find it convenient to formulate them for trigonometric polynomials, so set
\[
E_n[g]={\rm inf}_{{\deg}S\leq n}\,\|g-S\|_{[0,2\pi]},
\]
where the ${\rm inf}$ is taken over all trigonometric polynomials $S$ of degree at most $n$, $n\geq 1$, and $g:[0,2\pi]\to\mathbb R$ is a $2\pi$-periodic function.
The following holds. Let $0<\alpha<1$. Then, with constants independent of $g$ and $n$,
\[
E_n[g]=O(n^{-\alpha}),\ n\to\infty,\quad {\rm iff}\quad \omega(g,t)=O(t^{\alpha}),\ t\to 0^+,
\]
where $\omega(g,\cdot)$ is the modulus of continuity of $g$ on $[0,2\pi]$.
Thus the error of approximation of a $2\pi$-periodic function $g$ by trigonometric polynomials of degree at most $n\geq 1$ decays at rate $O(n^{-\alpha})$ iff $g$ satisfies a Lipschitz condition of order $0<\alpha<1$. More generally, if $k\geq 1$,
\[
E_n[g]=O(n^{-k-\alpha}),\ n\to\infty,\quad {\rm iff}\quad \omega\!\left(g^{(k)},t\right)=O(t^{\alpha}),\ t\to 0^+.
\]

\subsection{Restricted range inequalities}\label{sec:restricted}

The best way to motivate weighted Jackson and Bernstein theorems is to look carefully at the Freud weight, 
$w_{\lambda}(x):=\exp(-|x|^{\lambda})$ for $\lambda\geq 1$. Calculations show that the weighted 
polynomial $x^nw_{\lambda}(x)$ achieves its maximum modulus at the so-called Freud number $n^{\frac{1}{\lambda}}$. With this in mind, Freud and Nevai (see \cite{Mha}) showed that there are positive constants $C_1,C_2$ such that, uniformly for every polynomial $S$ of degree at most $n\geq 1$ and every $p\geq 1$,
\[
||Sw_{\lambda}||_{L_{p}\left(\mathbb R\right)}\leq C_1||Sw_{\lambda}||_{L_{p}\left(\left[-C_2n^{\frac{1}{\lambda}},C_2n^{\frac{1}{\lambda}}\right]\right)}.
\]
Outside, the interval $\left[-C_2n^{\frac{1}{\lambda}},
C_2n^{\frac{1}{\lambda}}\right]$, the function $fw_{\lambda}$ decays quickly to zero for large argument so one cannot hope to approximate
$fw_{\lambda}$ by $Sw_{\lambda}$. Indeed, a tail term, 
$||fw_{\lambda}||_{L_{p}\left(|x|\geq Cn^{\frac{1}{\lambda}}\right)}$
should typically appear in the weighted modulus of continuity for a weighted Jackson theorem. 
Ditzian and Lubinsky constructed a weighted modulus of continuity which has such a tail term (see \cite{Lub}, and the discussion below).
\medskip
Inequalities of the form 
\[
||Sw_{\lambda}||_{L_{p}\left(\mathbb R\right)}\leq C_1||Sw_{\lambda}||_{L_{p}\left(\left[-C_2n^{\frac{1}{\lambda}},C_2n^{\frac{1}{\lambda}}\right]\right)}.
\]
are called restricted range inequalities. The crucial observation is that the right hand side of the estimate involves a sequence of intervals which expand with $n$ as $n$ grows but are fixed for fixed $n$.
\medskip

Sharp forms of the restricted range inequality above exist for all $p\geq 1$ and for large classes of weights (see for example \cite{Mha,Lub1,ST}). For our purposes, we need the following standard result, due in this generality to Mhaskar, Rakhmanov and Saff; see \cite{Mha,Lub1,ST} for statements and proofs.
\medskip
Let $w=\exp(-\Phi):\mathbb R\to (0,\infty)$ with $\Phi$ even. We call $\Phi$ an external field for the weight $w$.
Assume $x\Phi'(x)$ is positive and increases on $(0,\infty)$ with limits $0$ and $\infty$ 
at $0$ and $\infty$. Then for each $n\geq 1$,
\[
\frac{2}{\pi}\int_{0}^{1}\frac{at\Phi'(at)}{\sqrt{1-t^2}}dt=n,\, a>0
\]
has a unique solution $a=a_n$, and for every polynomial $p_n$ of degree at most $n$, $n\geq 1$,
\[
\|p_nw\|_{\mathbb R}=\|p_nw\|_{[-a_n,a_n]}
\]
(see \cite{Mha,Lub1,ST}). The number $a_n$ is the Mhaskar--Rakhmanov--Saff number associated with the external field $\Phi$.
Consider again the even Freud weight $w_{\lambda}(x):=\exp(-|x|^{\lambda})$ for $\lambda\geq 1$. 

Then for every polynomial $p_n$ of degree at most $n\geq 1$
\[
\|p_nw_{\lambda}\|_{\mathbb R}=\|p_nw_{\lambda}\|_{[-a_n,a_n]}
\]
with $a=a_n:=(n\gamma_{\lambda})^{\frac{1}{\lambda}}$ and
$\gamma_{\lambda}:=\frac{\sqrt{\pi}\Gamma\left(\frac{\lambda}{2}\right)}{2\Gamma\left(\frac{\lambda+1}{2}\right)}$ with $\Gamma$ the Gamma function (Euler's integral of the second kind).  
\medskip

In particular:
\[
||p_n\exp(-|x|^2)||_{\mathbb R}=||p_n\exp(-|x|^2)||_{[-\sqrt{n},\sqrt{n}]}
\]
and
\[
||p_n\exp(-|x|)||_{\mathbb R}=||p_n\exp(-|x|)||_{[-\frac{\pi n}{2},\frac{\pi n}{2}]}.
\]

\subsection{Forward and Converse Jackson and Bernstein: Freud weights}

In formulating forward and converse Jackson and Bernstein theorems for Freud weights (see \cite{Lub}), we proceed as follows. 
For a real-valued function $f:[-1,1]\to \mathbb R$, we recall for $h>0$ and $k\geq 0$ the $k$th-order symmetric difference
\[
\Delta_{h}^{k}f(x)=\sum_{i=0}^k \binom{k}{i}(-1)^{i}\,f\!\left(x+\tfrac{kh}{2}-ih\right).
\]
(If any argument of $f$ lies outside $[-1,1]$, we set the difference to $0$.) 
\medskip

We use as our motivating example $w_{\lambda}(x):=\exp(-|x|^{\lambda})$ with $\lambda>1$ (see \cite{Lub}), and we restrict ourselves to $p=\infty$.
\medskip

The first-order weighted modulus of continuity is defined for suitable $f:\mathbb R\to \mathbb R$ as
\[
\omega_{1,\infty}(f,w_{\lambda},t):={\rm sup}_{0<h\leq t}\|w_{\lambda}\Delta_h(f)\|_{[-h^{\frac{1}{1-\lambda}},h^{\frac{1}{1-\lambda}}]}+{\rm inf}_{C\in \mathbb R}\|(f-C)w_{\lambda}\|_{\mathbb R\setminus ([-h^{\frac{1}{1-\lambda}},h^{\frac{1}{1-\lambda}}])}.
\]
The ${\rm inf}$ in the tail term ensures that if $f$ is constant the modulus vanishes identically, as one expects of a first-order modulus. Notice that substituting $h:=n^{-1+\frac{1}{\lambda}}$ gives $[-h^{\frac{1}{1-\lambda}},h^{\frac{1}{1-\lambda}}]=[-C_{1}a_n,C_{1}a_n]$.
\medskip

More generally, for $r\geq 1$, the $r$th-order weighted modulus of continuity for suitable $f:\mathbb R\to \mathbb R$ is
\[
\omega_{r,\infty}(f,w_{\lambda},t):={\rm sup}_{0<h\leq t}\|w_{\lambda}\Delta_{h}^r(f)\|_{[-h^{\frac{1}{1-\lambda}},h^{\frac{1}{1-\lambda}}]}+{\rm inf}_{{\rm deg}R\leq r-1}\|(f-R)w_{\lambda}\|_{\mathbb R\setminus([-h^{\frac{1}{1-\lambda}},h^{\frac{1}{1-\lambda}}])}.
\]

For $0<\alpha<r$, the following forward and converse estimates hold uniformly in $f$ and $n$ (see \cite{Lub}). First,
\[
E_n[f,w_{\lambda}]\leq C\,\omega_{r,\infty}\!\left(f,w_{\lambda},\tfrac{a_n}{n}\right)\leq C_1\left(\tfrac{a_n}{n}\right)^r\|f^{(r)}w_{\lambda}\|_{\mathbb R}.
\]
Second,
\[
\omega_{r,\infty}(f,w_{\lambda},t)=O\!\left(t^{\alpha}\right),\ t\to 0^+\quad {\rm iff}\quad
E_n[f,w_{\lambda}]=O\!\left(\left(\tfrac{a_n}{n}\right)^{\alpha}\right),\ n\to\infty.
\]
\medskip

The fractional order $\frac{a_n}{n}$ makes sense given it is the precise analogue of the fractional order $\frac{b-a}{n}$ in the classical case as above for a finite interval $[a,b]$.

Results such as the above (with different modulus) appear in \cite{Dam, Dam1, Dam2} for even weights of faster than smooth exponential decay for large argument and we refer the interested reader to these papers for more details. 

\section{Weighted deep polynomials}\label{sec:weighteddeep}

Moving forward, given a fixed weight $w$, a weighted deep polynomial is defined as
\[
Q(x) \;=\; \bigl[w(x)\bigr]^\alpha \cdot P(x),
\]
where the exponent \(\alpha\) is typically chosen to depend on \(D\) (for instance, \(\alpha = D\)). The rationale is that raising the weight to a high power localizes the effective region of approximation.
We note that as before, we will choose the weight $w$ to balance the growth of $P$ for large argument.

For different classes of real valued functions $f$ and weights $w$, we now show that 
numerically weighted deep polynomials capture $f$ extremely well in the uniform norm in particular, for large argument.
\medskip

We now study this problem in detail.

\section{Density and Existence of Best Approximants}\label{sec:existence}

Throughout this section $K\subset\mathbb R$ is a compact set with infinitely many points (e.g.\ a nondegenerate interval) and $w\colon K\to(0,\infty)$ is a fixed continuous, strictly positive weight. We measure error in the plain uniform norm $\|g\|_{\infty,K}:=\max_{x\in K}|g(x)|$, which is the norm used in our numerical experiments; the weight enters through the approximant, not the norm. For $D\in\mathbb N$ let
\[
\mathcal{Q}_D=\Bigl\{\,Q=[w]^{\alpha}\,(p_1\circ p_2\circ\cdots\circ p_L)\ :\ L\ge 1,\ \textstyle\prod_{\ell=1}^{L}\deg(p_\ell)\le D,\ \alpha\ge 0\,\Bigr\}
\]
denote the weighted deep polynomials of composite degree at most $D$.

\begin{theorem}[Density]\label{thm:WeightedDeepDensity}
For every $f\in C(K)$ and every $\varepsilon>0$ there exist $D\in\mathbb N$ and $Q\in\mathcal{Q}_D$ with $\|f-Q\|_{\infty,K}<\varepsilon$. That is, $\bigcup_{D\ge 1}\mathcal{Q}_D$ is dense in $C(K)$.
\end{theorem}

\begin{proof}
Fix $f\in C(K)$, $\varepsilon>0$, and any admissible exponent $\alpha\ge 0$. Since $w$ is continuous and strictly positive on the compact set $K$, the function $g:=f\,[w]^{-\alpha}$ is continuous on $K$. By the Weierstrass approximation theorem there is a polynomial $P$ with $\|g-P\|_{\infty,K}<\varepsilon/\|[w]^{\alpha}\|_{\infty,K}$. Taking $L=1$, $p_1=P$ (the identity outer layer), the weighted deep polynomial $Q=[w]^{\alpha}P$ lies in $\mathcal{Q}_D$ with $D=\deg P$ and satisfies
\[
\|f-Q\|_{\infty,K}=\bigl\|[w]^{\alpha}(g-P)\bigr\|_{\infty,K}\le \|[w]^{\alpha}\|_{\infty,K}\,\|g-P\|_{\infty,K}<\varepsilon. \qedhere
\]
\end{proof}

\begin{theorem}[Existence of a best approximant]\label{thm:WeightedDeepExist}
Fix $D\in\mathbb N$, a finite list of layer-degree patterns $(d_1,\dots,d_L)$ with $\prod_\ell d_\ell\le D$, and $\alpha\ge0$, and let $\mathcal{Q}_D'\subseteq\mathcal{Q}_D$ be the corresponding family. For every $f\in C(K)$ the value
\[
E^{*}:=\inf_{Q\in\mathcal{Q}_D'}\|f-Q\|_{\infty,K}
\]
is attained by some $Q^{*}=[w]^{\alpha}P^{*}$ with $\deg P^{*}\le D$, lying in the uniform closure of $\mathcal{Q}_D'$ in $C(K)$.
\end{theorem}

\begin{proof}
Let $Q_k=[w]^{\alpha}P_k\in\mathcal{Q}_D'$ satisfy $\|f-Q_k\|_{\infty,K}\to E^{*}$. Then $(\|Q_k\|_{\infty,K})_k$ is bounded, and since $\min_K[w]^{\alpha}>0$, so is $(\|P_k\|_{\infty,K})_k$. Each $P_k$ lies in the space $\mathcal P_D$ of polynomials of degree at most $D$; because $K$ has infinitely many points, $\|\cdot\|_{\infty,K}$ is a norm on the finite-dimensional space $\mathcal P_D$, so a subsequence of $(P_k)$ converges uniformly on $K$ to some $P^{*}\in\mathcal P_D$. The limit $Q^{*}:=[w]^{\alpha}P^{*}$ satisfies $\|f-Q^{*}\|_{\infty,K}=E^{*}$.
\end{proof}

\begin{remark}[Non-closedness of a fixed architecture]\label{rem:nonclosed}
The infimum need not be attained inside $\mathcal{Q}_D'$ itself, because the set of compositions of a fixed pattern is not closed in $C(K)$. For the pattern $(d_1,1)$ with $d_1\ge2$, the normalized compositions $q_k\circ p_k$ with $p_k(x)=x^{d_1}+kx$ (monic, zero constant term) and $q_k(y)=y/k$ converge uniformly on $K$ to the identity, which is not a composition of that pattern. Such degeneration of layers is the standard obstruction to existence in nonlinear approximation \cite{Braess}; in computation it appears as unbounded drift of layer coefficients. The monic normalization limits this drift in practice, and the fine-tuning construction of Section~\ref{sec:finetune} avoids the issue entirely, since there the inner layers are fixed and the outer fit is a convex problem whose minimum is attained.
\end{remark}

\begin{remark}[On uniqueness]
Best uniform approximants from the nonlinear family $\mathcal{Q}_D'$ need not be unique: the classical Haar-system uniqueness theory applies to linear spaces only. Uniqueness can be restored in strictly convex norms (for instance $L^{p}$ with $1<p<\infty$); since our algorithms use random restarts and retain the best minimizer found, non-uniqueness poses no practical difficulty.
\end{remark}

\section{Asymmetric weighted approximation}\label{sec:asymmetric}

\subsection{Reduction to a compact interval}\label{subsec:compact}

The role of the one-sided weight can be stated precisely: it converts uniform approximation on a half line into uniform approximation on a compact interval whose length grows slowly with the degree, in direct analogy with the restricted-range inequalities of Section~\ref{sec:restricted}. We record an elementary, self-contained version adapted to one-sided fields; it requires nothing of the weight beyond an upper bound on the decaying side.

\begin{lemma}\label{lem:chebgrowth}
Let $P$ be a polynomial of degree at most $n$ and let $a<s$. For $x\ge s$,
\[
|P(x)|\;\le\;\|P\|_{[a,s]}\,\Bigl(\xi(x)+\sqrt{\xi(x)^{2}-1}\Bigr)^{n},
\qquad \xi(x)=\frac{2x-a-s}{s-a}\,.
\]
\end{lemma}

\begin{proof}
This is the classical extremal growth property of the Chebyshev polynomial: among all polynomials of degree at most $n$ bounded by $1$ on $[a,s]$, $|P(x)|\le|T_n(\xi(x))|$ for $x\notin[a,s]$, and $T_n(\xi)=\cosh(n\,\operatorname{arccosh}\xi)\le\bigl(\xi+\sqrt{\xi^{2}-1}\bigr)^{n}$ for $\xi\ge1$; see \cite{T}.
\end{proof}

\begin{proposition}[Compactification]\label{prop:compact}
Let $a<s$, $c>0$, $\beta>1$, and let $w:[a,\infty)\to(0,1]$ be continuous with
\[
w(x)\;\le\;e^{-c\,(x-s)^{\beta}}\qquad\text{for }x\ge s.
\]
There exist $K=K(a,s)>0$ and $n_0=n_0(a,s,c,\beta)$ such that for all $n\ge n_0$, setting
\[
b_n \;=\; s+\Bigl(\frac{K\,n\log n}{c}\Bigr)^{1/\beta},
\]
every polynomial $P$ of degree at most $n$ satisfies
\[
\|wP\|_{[b_n,\infty)}\;\le\;e^{-n}\,\|P\|_{[a,s]}\,.
\]
Consequently, for every $f\in C[a,\infty)$,
\[
\|f-wP\|_{[a,b_n]}\;\le\;\|f-wP\|_{[a,\infty)}\;\le\;\|f-wP\|_{[a,b_n]}+\|f\|_{[b_n,\infty)}+e^{-n}\,\|P\|_{[a,s]}\,.
\]
\end{proposition}

\begin{proof}
For $x\ge s$ write $v=x-s\ge0$. By Lemma~\ref{lem:chebgrowth} and $\xi+\sqrt{\xi^{2}-1}\le2\xi$,
\[
w(x)\,|P(x)|\;\le\;\|P\|_{[a,s]}\,\exp\bigl(n\log(2\xi(x))-c\,v^{\beta}\bigr),
\]
and $\xi(x)=1+2v/(s-a)\le C_{0}(1+v)$ with $C_{0}=1+2/(s-a)$. Set $\psi(v)=n\log\bigl(2C_{0}(1+v)\bigr)-c\,v^{\beta}$. Then
\[
\psi'(v)=\frac{n}{1+v}-c\beta\,v^{\beta-1}<0
\qquad\text{whenever}\qquad c\beta\,v^{\beta-1}(1+v)>n,
\]
in particular for all $v\ge\bigl(n/(c\beta)\bigr)^{1/\beta}$. Let $v_n=\bigl(2n\log(2C_{0}(1+n))/c\bigr)^{1/\beta}$. Since $\beta>1$ we have $v_n\le n$ for all $n$ large, hence
\[
\psi(v_n)\;\le\;n\log\bigl(2C_{0}(1+n)\bigr)-2n\log\bigl(2C_{0}(1+n)\bigr)\;=\;-\,n\log\bigl(2C_{0}(1+n)\bigr)\;\le\;-n,
\]
and $v_n\ge(n/(c\beta))^{1/\beta}$ for $n$ large, so $\psi$ is decreasing on $[v_n,\infty)$. Therefore $\sup_{x\ge s+v_n}w(x)|P(x)|\le e^{-n}\|P\|_{[a,s]}$. Since $2\log(2C_{0}(1+n))\le K\log n$ for $K=K(a,s)$ and $n$ large, the choice $b_n\ge s+v_n$ of the statement inherits the bound by monotonicity of $\psi$. The two-sided error estimate follows from the triangle inequality $\|f-wP\|_{[b_n,\infty)}\le\|f\|_{[b_n,\infty)}+\|wP\|_{[b_n,\infty)}$ and monotonicity of the supremum norm in the domain.
\end{proof}

\begin{remark}\label{rem:compact}
(i) On the decaying side the proposition plays the role of the restricted-range identity $\|p_nw\|_{\mathbb R}=\|p_nw\|_{[-a_n,a_n]}$ of Section~\ref{sec:restricted}, with $b_n$ at the Mhaskar--Rakhmanov--Saff scale $n^{1/\beta}$ up to a logarithmic factor. Sharper forms, with exponentially small tails beyond $(1+\delta)a_n$ and no logarithm, are available from the theory cited there \cite{Mha,Lub1,ST}; the elementary bound above suffices for our purposes and covers one-sided fields directly.
(ii) If in addition $w\equiv1$ on $[a,s]$, then $\|P\|_{[a,s]}=\|wP\|_{[a,s]}$ and the tail is controlled by the weighted approximant itself. The smooth one-sided weight $w(x)=\exp\bigl(-c[\log(1+e^{x-s})]^{\beta}\bigr)$ of Section~\ref{sec:finetune} satisfies the decay hypothesis, since $\log(1+e^{t})\ge t$ for $t\ge0$, and is bounded below by the constant $e^{-c(\log2)^{\beta}}$ on $(-\infty,s]$.
(iii) The practical consequence is that once the working interval contains $[a,b_n]$, its right endpoint is immaterial: beyond $b_n$ the error of a weighted approximant is dominated by the tail of $f$ itself. This is confirmed numerically in Section~\ref{sec:finetune}, where the weighted error falls below machine precision past the transition region.
\end{remark}

\subsection{A model problem}\label{subsec:model}

Consider approximating the function
\[
f(x) = e^{-x},
\]
which exhibits contrasting behavior on the two halves of \(\mathbb{R}\): as \(x\to-\infty\) the function diverges to \(+\infty\) (since \(e^{-x}=e^{|x|}\)), while for \(x\ge0\) it decays to \(0\). To obtain a finite approximation error over a domain that extends far to the right, the approximant must mimic this behavior.

We define the weight via a ReLU-based function:
\[
w\bigl(\operatorname{ReLU}(x)\bigr)=
\begin{cases}
1, & x<0,\\[1mm]
e^{-x^2}, & x\ge0.
\end{cases}
\]
Then, the weighted deep polynomial approximant is given by
\[
Q_n(x) \;=\; \Bigl[w\bigl(\operatorname{ReLU}(x)\bigr)\Bigr]^n\,P(x),
\]
where the exponent \(n\) is tied to the degree of the deep polynomial. For \(x<0\) the weight is constant so that \(Q_n(x)=P(x)\) can grow to match the divergence of \(e^{-x}\); for \(x\ge0\) the weight becomes
\[
e^{-n x^2},
\]
which decays rapidly as \(n\) increases. In effect, significant contributions for \(x\ge0\) occur only when \(n\,x^2\) is small, i.e.\ for \(x\) of order \(1/\sqrt{n}\).

A change of variables \(y=\sqrt{n}\,x\) (so that \(dx=dy/\sqrt{n}\)) transforms the error on \([0,\infty)\) to
\[
E_n^+ = \frac{1}{\sqrt{n}}\int_{0}^{\infty} \Bigl|e^{-y/\sqrt{n}}-P\Bigl(\frac{y}{\sqrt{n}}\Bigr)e^{-y^2}\Bigr|^2\,dy.
\]
Since \(e^{-y/\sqrt{n}}\) converges uniformly to 1 on any fixed interval as \(n\to\infty\), the dominant contribution arises from a bounded region in \(y\). Appropriately choosing \(P(x)\) (via rescaling) thus ensures that the error in approximating \(f(x)\) can be made arbitrarily small.

\subsection{Effective Interval and Endpoint Localization}

Due to the rapid decay of the weight for \(x\ge0\), the approximation is effectively localized to a small interval near the origin. After rescaling with \(y=\sqrt{n}\,x\), the significant region is confined to an interval \([0,a]\) in the \(y\)-variable. This endpoint \(a\) is determined by the condition
\begin{equation}\label{eq:endpoint}
\int_0^a \frac{\Phi'(t)}{\sqrt{a^2-t^2}}\,dt \;=\; \frac{\pi}{2},
\end{equation}
where \(\Phi(t)=t^2\) plays the role of an external field. Changing variables via \(t=a u\) yields
\[
a\int_0^1 \frac{2a u}{a\sqrt{1-u^2}}\,du
=\; 2a\int_0^1 \frac{u}{\sqrt{1-u^2}}\,du
=\; 2a.
\]
Thus, condition \eqref{eq:endpoint} implies
\[
2a \;=\; \frac{\pi}{2} \quad\Longrightarrow\quad a=\frac{\pi}{4}.
\]
In the original \(x\)-variable, the effective interval is then approximately
\[
\left[0,\frac{\pi}{4\sqrt{n}}\right].
\]
This localization is key to ensuring that the weighted deep polynomial \(Q_n(x)\) converges to \(f(x)\) as \(n\to\infty\).

\section{Computational Results}\label{sec:computational}

\subsection{Algorithm for Deep Weighted Polynomial Approximation}\label{subsec:deep-weighted-poly}

We seek to approximate a real-valued function \(f(x)\) by a \emph{deep weighted polynomial} of the form
\[
Q(x) \;=\; w(x)^\gamma\,h^{(L)}(x),
\]
where the weight \(w(x)\colon[a,b] \to \mathbb{R}_{>0}\) is typically chosen to enforce exponential decay or localization, and \(\gamma \ge 0\) is a tunable exponent. The function \(h^{(L)}(x)\) is the output of a scalar-valued \emph{computational graph} of depth \(L\), where each layer is constructed recursively via the polynomial update
\[
h^{(\ell)}(x) \;=\; \sum_{i=0}^{m_\ell - 1} a^{(\ell)}_i \left(h^{(\ell-1)}(x)\right)^i, \quad h^{(0)}(x) = x,
\]
and the parameter vector \(\theta = (a^{(1)}, a^{(2)}, \dots, a^{(L)}) \in \mathbb{R}^{n_{\mathrm{deep}}}\) comprises all coefficients. Each inner layer ($\ell < L$) is normalized to be monic with zero constant term, exactly as in Section~\ref{sec:intro}, so the effective parameter count is
\[
n_{\mathrm{deep}} \;=\; \Bigl(\sum_{\ell=1}^L d_\ell\Bigr) - L + 2,
\]
where $d_\ell = m_\ell - 1$ is the degree of layer $\ell$ and $m_\ell$ its width; this is the count used to match budgets against classical baselines. The construction generalizes the computational-graph framework of Jarlebring et al.~\cite{Jar} from matrix functions to nonlinear polynomial approximants over one-dimensional input, optimizing $\theta$ to minimize the weighted least-squares loss
\[
L(\theta) \;=\; \sum_{i=1}^N \left(Q(x_i; \theta) - f(x_i)\right)^2 \Delta x,
\]
where \(\{x_i\}_{i=1}^N \subset [a,b]\) is a uniform sampling grid and \(\Delta x = (b - a)/N\) is the quadrature weight. The loss is minimized using gradient-based optimization, and gradients \(\nabla_\theta L\) are computed efficiently by automatic differentiation through the computational graph.

\smallskip

Algorithm~\ref{alg:deep-weighted-general} summarizes the procedure. A scalar DAG is constructed with input node \(x = h^{(0)}\), followed by composition layers of prescribed widths \(\{m_\ell\}\), with inner layers normalized as above. The graph consists only of multiplication and linear-combination operations, making it amenable to automatic differentiation. To mitigate local minima, we perform \(R\) random restarts with Gaussian initialization; for each restart the loss is minimized by backpropagating residuals through the graph \cite{Jar}, and the best coefficient vector \(\theta^*\) is retained. The weight \(w(x)^\gamma\) focuses accuracy in regions where \(w(x)\) is large without increasing global polynomial complexity.

\paragraph{Cost model.} Two budgets are relevant when comparing approximants of different structure: the number of \emph{trainable parameters}, i.e.\ coefficients adjusted by the fitting procedure, and the \emph{evaluation cost}, i.e.\ the number of multiplications required to evaluate the approximant once. For a single polynomial of degree $n$ the two essentially coincide: $n+1$ parameters and $n$ multiplications by Horner's rule (or Clenshaw's recurrence in Chebyshev form). For composite approximants they diverge. A fixed map, such as the Newton update for the sign function in Section~\ref{sec:intro}, has no trainable parameters yet a nonzero per-evaluation cost, while a deep polynomial with layer degrees $(d_1,\dots,d_L)$ evaluates in $O\bigl(\sum_\ell d_\ell\bigr)$ multiplications but realizes effective degree $\prod_\ell d_\ell$. Reporting trainable parameters alone would therefore flatter methods with fixed structure. Throughout the paper we match trainable parameters across methods and, in Section~\ref{sec:finetune} where part of the structure is fixed, we verify in addition that evaluation costs are of the same order.

\begin{algorithm}[h]
\caption{Deep Weighted Polynomial via Computational Graphs}
\label{alg:deep-weighted-general}
\DontPrintSemicolon
\KwIn{$f\colon[a,b]\to\mathbb R$, weight $w\colon[a,b]\to\mathbb R_{>0}$, exponent $\gamma \ge 0$, samples $\{x_i\}_{i=1}^N$, layers $L$, widths $\{m_\ell\}$, restarts $R$, step $\eta$}
\KwOut{$\theta^* = (a^{(1)}, \dots, a^{(L)}) \in \mathbb{R}^{n_{\mathrm{deep}}}$, where $n_{\mathrm{deep}}=(\sum_\ell d_\ell)-L+2$ (normalized count)}
Let $n_{\mathrm{deep}} = \bigl(\sum_{\ell=1}^L d_\ell\bigr) - L + 2$ (inner layers normalized: monic, zero constant term)\;

\textbf{Graph construction:} Create DAG $\mathcal{G}$ with input $h^{(0)} = x$\;
\For{$\ell = 1$ \KwTo $L$}{
  Add nodes $(h^{(\ell-1)})^i$, $i=0,\dots,m_\ell - 1$\;
  Add one linear-combination node: $h^{(\ell)} = \sum_{i=0}^{m_\ell - 1} a^{(\ell)}_i (h^{(\ell-1)})^i$\;
}
Final output node computes $Q(x) = w(x)^\gamma\, h^{(L)}(x)$\;

Set $L^* \gets +\infty$\;

\For{$r = 1$ \KwTo $R$}{
  Initialize $\theta \sim \mathcal{N}(0, I)$\;
  
  \Repeat{convergence}{
    Compute loss: $L(\theta) = \sum_i \left[ Q(x_i; \theta) - f(x_i) \right]^2 \, \Delta x$\;
    Compute $\nabla_\theta L$ by backpropagation through $\mathcal{G}$~\cite{Jar}\;
    Update parameters: $\theta \gets \theta - \eta \nabla_\theta L$\;
  }
  \If{$L(\theta) < L^*$}{
    $\theta^* \gets \theta$, $L^* \gets L(\theta)$\;
  }
}
\KwRet{$\theta^*$}
\end{algorithm}

\medskip

\noindent\textbf{Implementation details.}  In our experiments we take \(N\approx600\) uniformly 
spaced sample points on \([a,b]\), set \(\gamma=1\) for the weighted variant and \(\gamma=0\) for
the unweighted baseline, and choose \(R=5\) random restarts.  Convergence is declared when the
relative change in \(L\) drops below \(10^{-12}\).  This framework readily extends to other weight
functions and to different loss functions (e.g.\ \(L_\infty\) or relative error), and forms the basis
for our deep weighted polynomial approximations throughout the paper.

\subsection{Exponential Function}

Polynomials fail to approximate functions on unbounded domains because any non-constant polynomial \( p(x) \) diverges to \(\pm\infty\) as \(x \to \pm\infty\). For example, only the constant polynomial can approximate \(e^x\) with finite error on \((-\infty, 0]\); therefore, the minimax error is no less than \(1/2\), since \(e^0 = 1\) while \(e^{-\infty} = 0\). This challenge highlights the limitations of classical approximation theory in settings where the target function exhibits extreme or asymmetric behavior at infinity.

As a first example, we consider a function that exhibits two distinct behaviors: it diverges exponentially as \(x \to -\infty\) and decays rapidly as \(x \to +\infty\). The target \(e^{-x}\) is itself entire (smooth, with all derivatives continuous); the non-smoothness that the deep polynomial must resolve is \emph{not} a property of \(f\) but is introduced by the one-sided weight described below, which glues a constant branch on \(x<0\) to a decaying branch on \(x\ge0\). It is precisely this non-smooth interface, together with the rapid right-tail decay, that makes the construction well-suited to deep polynomial approximants.

In our approach, the weight function is defined asymmetrically as
\[
w(x) =
\begin{cases}
1, & x < 0, \\[1mm]
e^{-x^2}, & x \ge 0,
\end{cases}
\]
so that the overall approximant takes the form
\[
Q(x) = w(x) \cdot q\bigl(p(x)\bigr).
\]
The asymmetry of $w$ is deliberate and is the whole point of the construction: we keep $w\equiv 1$ on $x<0$ so that the composite polynomial is free to reproduce the growth $e^{|x|}$ there, while we let $w$ decay on $x\ge0$ to tame the right tail. One could of course replace this $C^1$ gluing by a globally $C^\infty$ one-sided weight (for instance a smooth sigmoidal transition), and nothing in our method requires the interface to be non-smooth; we use the simplest ReLU-type gluing because it is convenient and because it lets us exhibit, in a controlled way, that the deep polynomial composition is robust to a non-smooth interface. Concretely, although \(w\) is continuous with \(w'(0^-) = w'(0^+) = 0\), its second derivative jumps at \(x=0\):
\[
w''(x) =
\begin{cases}
0, & x < 0, \\
(4x^2 - 2) e^{-x^2}, & x \ge 0,
\end{cases}
\]
so that \(w''(0^+) = -2\neq 0=w''(0^-)\). The composition \(q(p(x))\) must therefore adapt across this interface; we report below that it does so without difficulty.

To ensure that the deep polynomial \(q(p(x))\) accurately tracks \(f(x) = e^{-x}\) across both regions, it must satisfy
\[
Q(x) =
\begin{cases}
q(p(x)) \approx e^{-x}, & x < 0, \\[1mm]
e^{-x^2} \cdot q(p(x)) \approx e^{-x}, & x \ge 0,
\end{cases}
\]
which is equivalent to enforcing
\[
q(p(x)) \approx \frac{e^{-x}}{e^{-x^2}} = e^{-x + x^2}
\qquad \text{for } x \ge 0.
\]
Thus, the same deep polynomial \(q(p(x))\) must serve the dual role of approximating \(e^{-x}\) on the left and \(e^{-x+x^2}\) on the right. This type of behavior cannot be captured by a single unweighted polynomial of modest degree, given minmax approximation relies on equioscillation, but arises naturally from the compositional structure of our deep graph-based approximation. The ability to resolve such non-smoothness is a key strength of the deep polynomial approach \cite{KY}. 

\begin{figure}[H]
    \centering
    \begin{minipage}{0.50\textwidth}
        \centering
        \includegraphics[width=\linewidth]{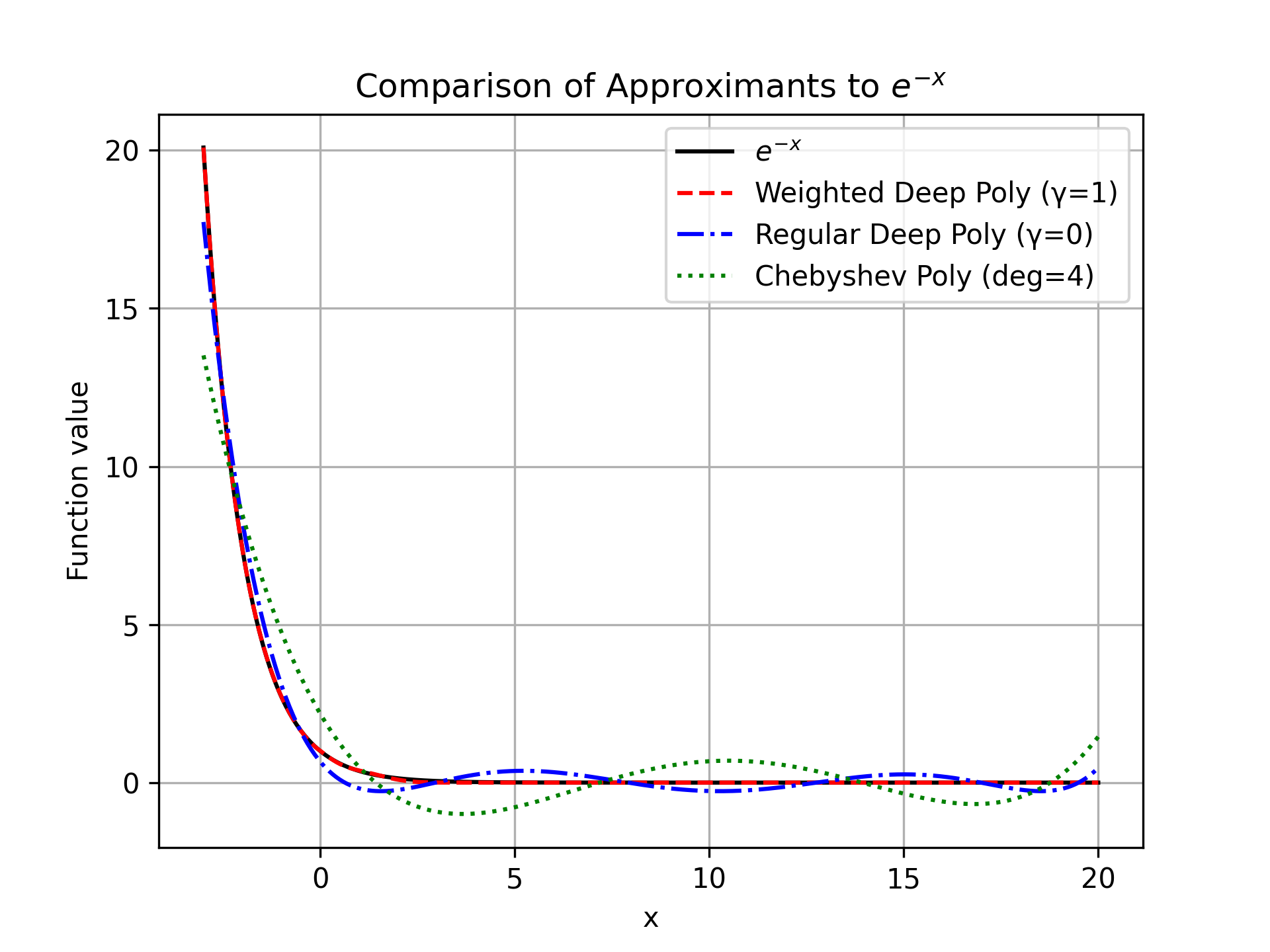}
        \caption*{(a) Exponential function with one-sided Gaussian weight}
    \end{minipage}\hfill
    \begin{minipage}{0.50\textwidth}
        \centering
        \includegraphics[width=\linewidth]{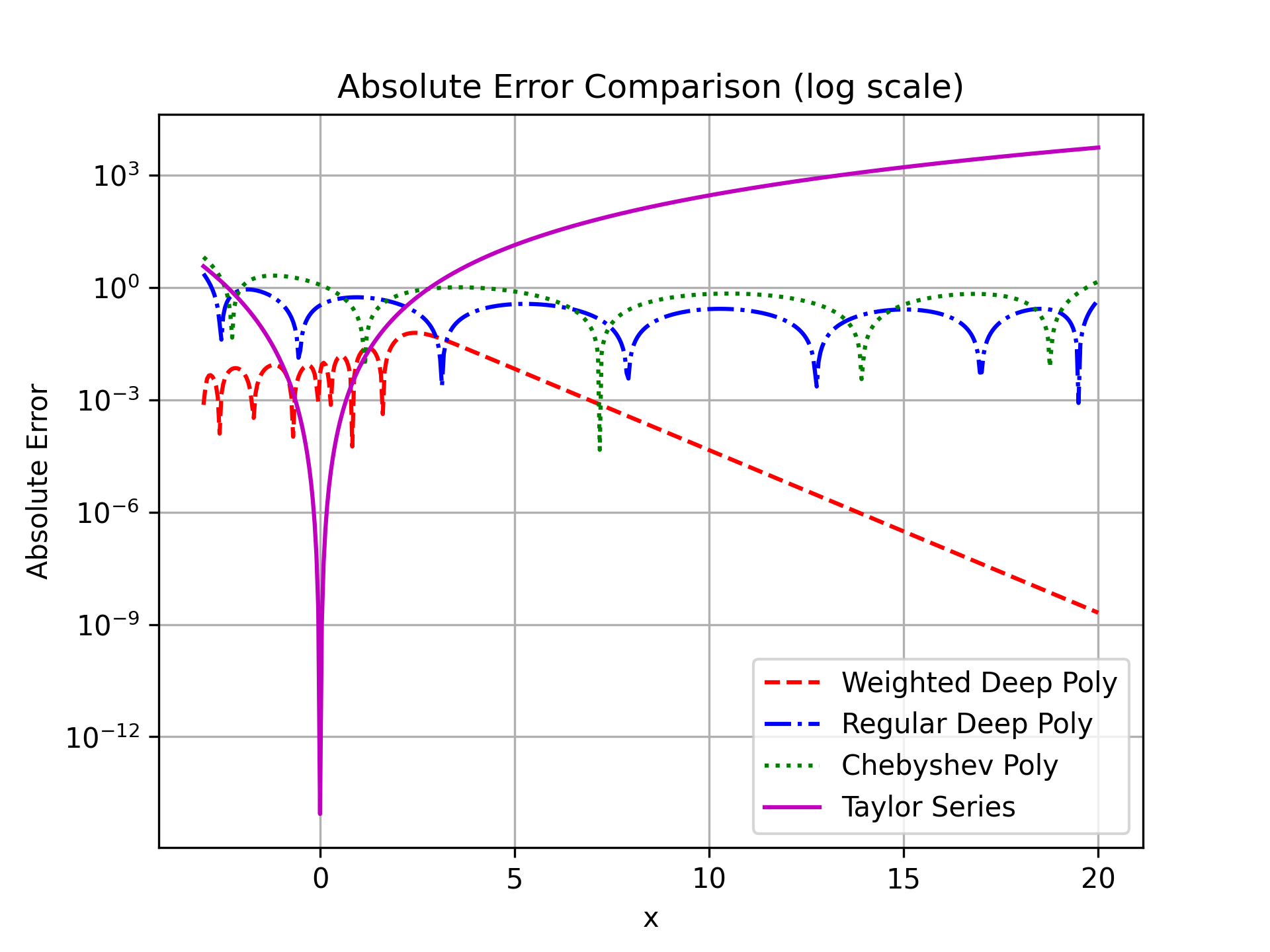}
        \caption*{(b) Pointwise absolute Error (log scale)}
    \end{minipage}
    \caption{Approximants for $e^{-x}$ on $[-5,20]$ at matched parameter budgets: the weighted deep polynomial against a Chebyshev polynomial baseline and an unweighted deep polynomial, all with $n_{\mathrm{deep}}$ free parameters as defined in Section~\ref{sec:computational}. Panel (b): pointwise absolute error on a logarithmic axis.}
    \label{fig:exponential}
\end{figure}

Figure~\ref{fig:exponential} compares the approximants for \(e^{-x}\) on \([-5, 20]\) at equal parameter budgets. The advantage of the weighted deep polynomial is most pronounced on the decaying side ($x \ge 0$), where the weight is active and where an unweighted polynomial of any fixed budget must trade accuracy on the tail against accuracy on the growing side; this trade-off is exactly the mechanism quantified by Proposition~\ref{prop:compact}.

\begin{figure}[H]
    \centering
    \begin{minipage}{0.45\textwidth}
        \centering
        \includegraphics[width=\linewidth]{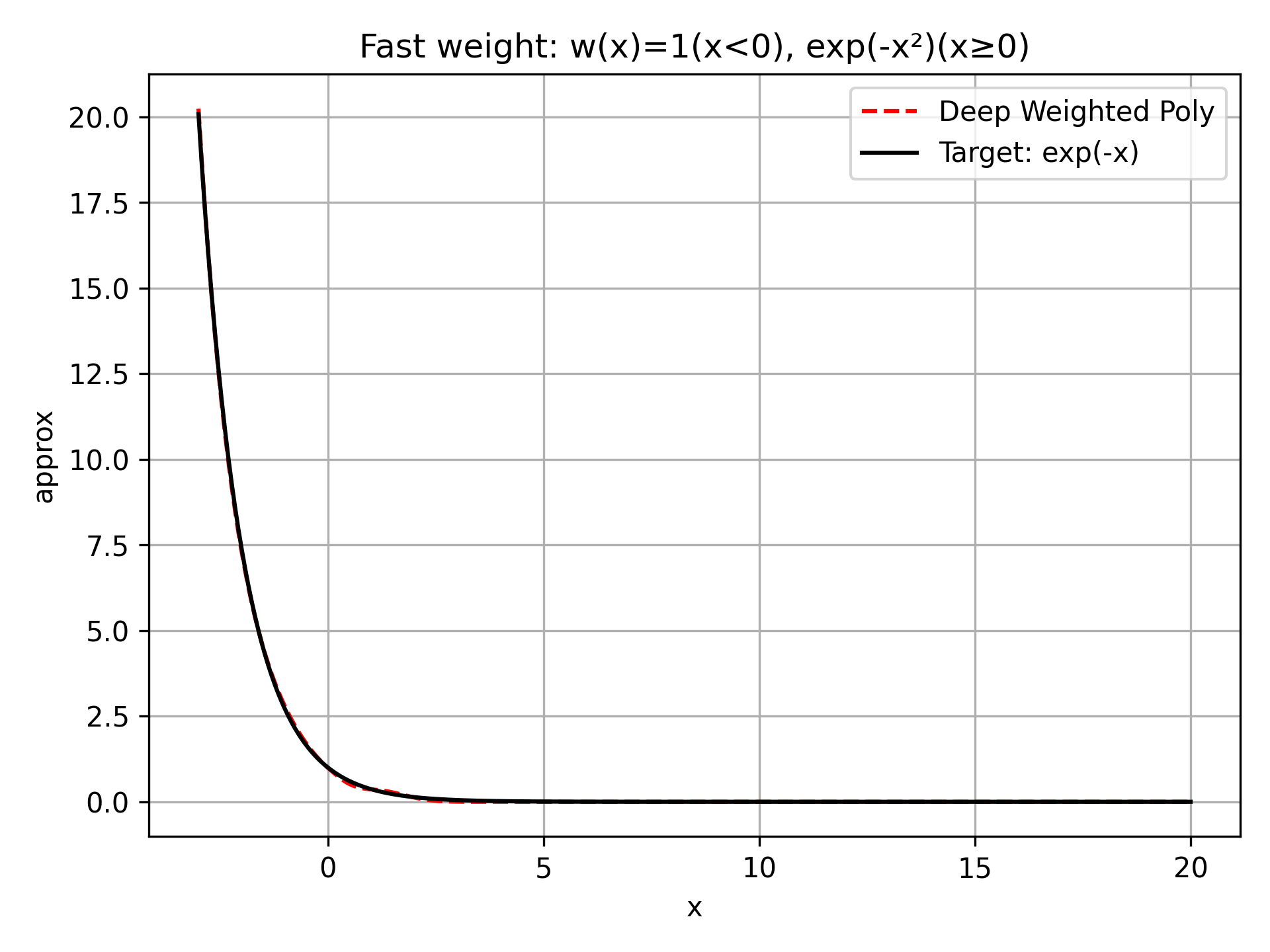}
        \caption*{(a) Deep weighted polynomial approximant for \(e^{-x}\) on \([-3,20]\) with fast weight \(w(x)=\begin{cases}
            1, & x<0,\\
            e^{-x^2}, & x\ge0.
        \end{cases}\)}
    \end{minipage}\hfill
    \begin{minipage}{0.45\textwidth}
        \centering
        \includegraphics[width=\linewidth]{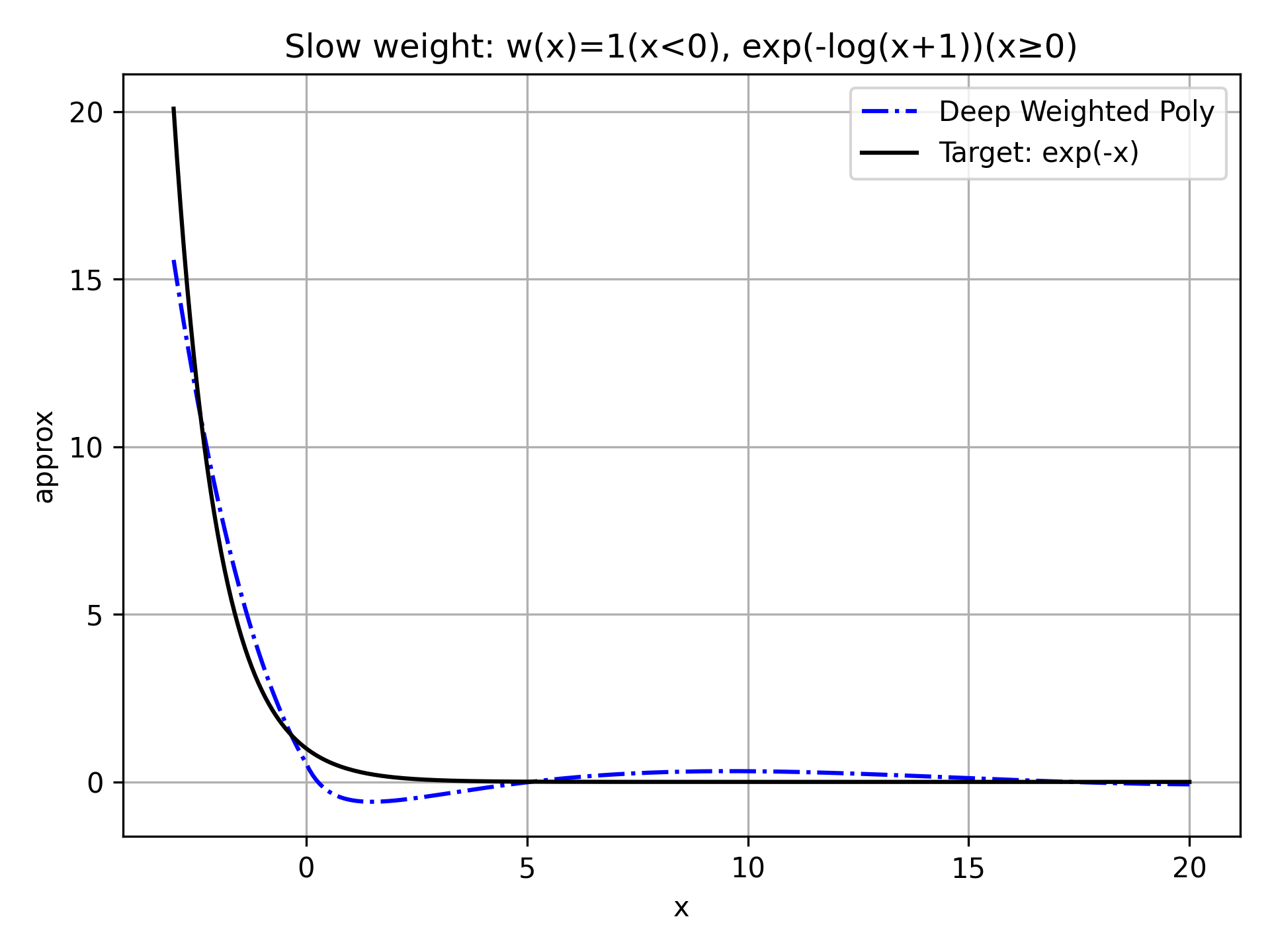}
        \caption*{(b) Deep weighted polynomial approximant for \(e^{-x}\) on \([-3,20]\) with slow weight \(w(x)=\begin{cases}
            1, & x<0,\\
            \tfrac{1}{x+1}, & x\ge0.
        \end{cases}\)}
    \end{minipage}
    \caption{Deep weighted polynomial approximants for \(e^{-x}\) on \([-3,20]\) with a one-sided Gaussian weight (left) and a one-sided reciprocal weight (right), against baselines of the same parameter budget.}
    \label{fig:weight_comparison}
\end{figure}

Figure~\ref{fig:weight_comparison} illustrates the crucial role played by the one-sided weight in the construction of deep polynomial approximants. If the weight decays too rapidly, the approximant receives little information from the right tail of the function, while if it decays too slowly, the polynomial must work harder to cancel the residual, leading to poor approximation. The Gaussian weight strikes an effective balance by suppressing the tail without overly penalizing the polynomial. We will later demonstrate in Section \ref{sec:generalizedweight} that this choice is nearly optimal.

\subsection{Airy Function}

We consider the Airy function of the second kind, denoted \(\mathrm{Bi}(x)\), which arises in a variety of physical applications. It satisfies the second-order linear differential equation
\[
y''(x) - x y(x) = 0,
\]
and admits the integral representation
\[
\mathrm{Bi}(x) = \frac{1}{\pi} \int_0^\infty \left[ \exp\left( \frac{1}{3} t^3 + x t \right) + \sin\left( \frac{1}{3} t^3 + x t \right) \right]\,dt.
\]
The function \(\mathrm{Bi}(x)\) grows rapidly for large positive \(x\), while exhibiting oscillatory behavior for negative \(x\). This asymmetry and analytic complexity make it a compelling test case for structured approximants.

Airy functions appear prominently in semiclassical quantum mechanics, especially in potential well and tunneling problems, and also in wave propagation near turning points. Due to their special-function status, they are generally computed using library routines or special-purpose quadrature, making them computationally expensive to evaluate repeatedly in applications.

We construct approximants of the form
\[
Q(x) = w(x)^\gamma \, q(p(x)),
\]
where \(p(x)\) is a normalized inner polynomial, \(q(\cdot)\) is an outer polynomial, and \(w(x)\) is a weight selected to suppress growth or oscillation in regions of lesser importance. There is no contradiction between using a fixed coefficient budget and claiming efficiency: by adapting the composite structure to the behavior of \(\mathrm{Bi}(x)\), the deep weighted polynomial realizes an \emph{effective (composite) degree} of $\prod_\ell(m_\ell-1)$ that far exceeds the degree of a single classical polynomial with the same number of free coefficients. It is in this sense, few free coefficients but high effective degree, that the construction is parameter-efficient, and all approximants compared in a given figure share the same raw coefficient count.

\begin{figure}[H]
    \centering
    \begin{minipage}{0.49\textwidth}
        \centering
        \includegraphics[width=\linewidth]{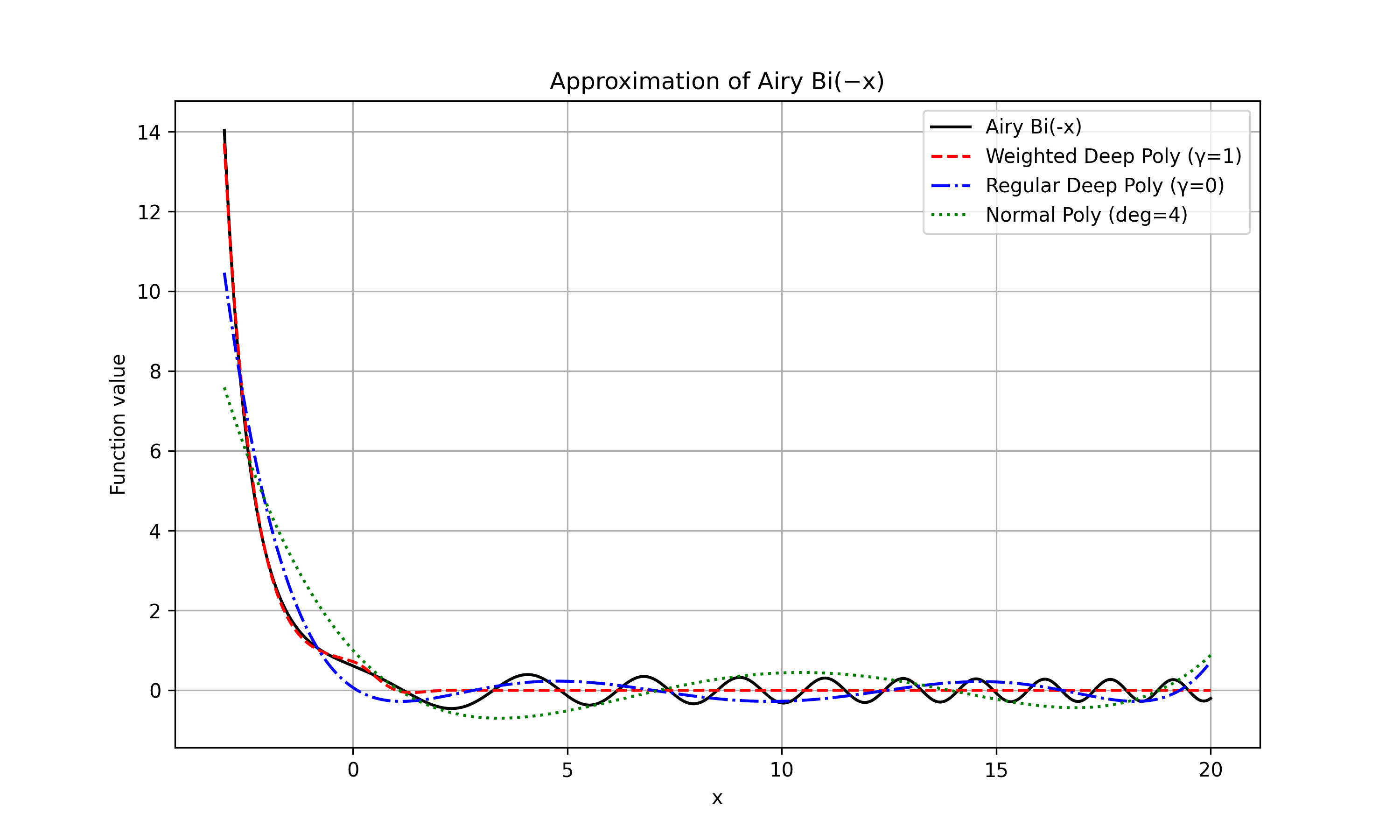}
        \caption*{(a) Airy function with one-sided Gaussian weight}
    \end{minipage}\hfill
    \begin{minipage}{0.49\textwidth}
        \centering
        \includegraphics[width=\linewidth]{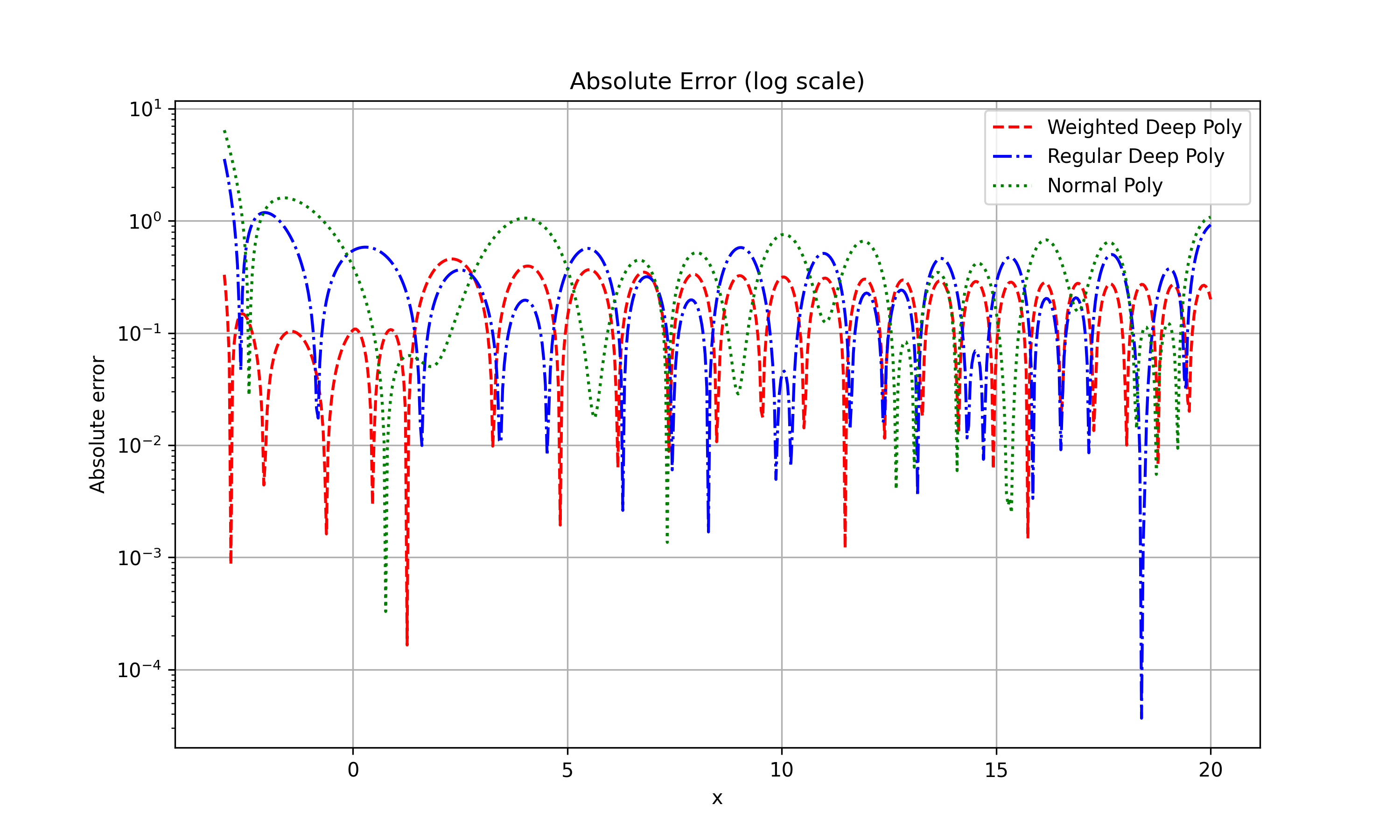}
        \caption*{(b) Pointwise absolute Error (log scale)}
    \end{minipage}
    \caption{Approximants for $\mathrm{Bi}(-x)$ on $[-3,20]$ at matched parameter budgets: the weighted deep polynomial, the unweighted deep polynomial, and the Chebyshev polynomial baseline, all with $n_{\mathrm{deep}}$ free parameters. Panel (b): pointwise absolute error on a logarithmic axis.}
    \label{fig:airy_approximants}
\end{figure}

Figure~\ref{fig:airy_approximants} shows the approximants, all with $n_{\mathrm{deep}}$ free parameters, for \(\mathrm{Bi}(-x)\) on \([-3, 20]\). Among the methods shown, the weighted deep polynomial best captures the growth across the domain. At this very small budget all methods, including the weighted deep polynomial, lose accuracy in the far right tail; the weighted approximant degrades most slowly. Recovering uniform tail accuracy requires a larger budget or a better-matched weight, which is addressed by the weight learning of Section~\ref{sec:generalizedweight} and, decisively, by the fine-tuned construction of Section~\ref{sec:finetune}, which resolves the entire decaying tail to machine precision on harder targets.

\subsection{Generalized symmetric Freud weight and weight learning}\label{sec:generalizedweight}

The restricted-range identity recalled in Section~\ref{sec:restricted} (the existence of a unique Mhaskar--Rakhmanov--Saff number $a_n$ with $\|p_nw\|_{\mathbb R}=\|p_nw\|_{[-a_n,a_n]}$) holds for any even external field satisfying the hypotheses stated there; we do not repeat it. Here we use it to motivate \emph{learning} the weight rather than fixing it in advance.

\begin{figure}[H]
    \centering
    \includegraphics[width=0.4\textwidth]{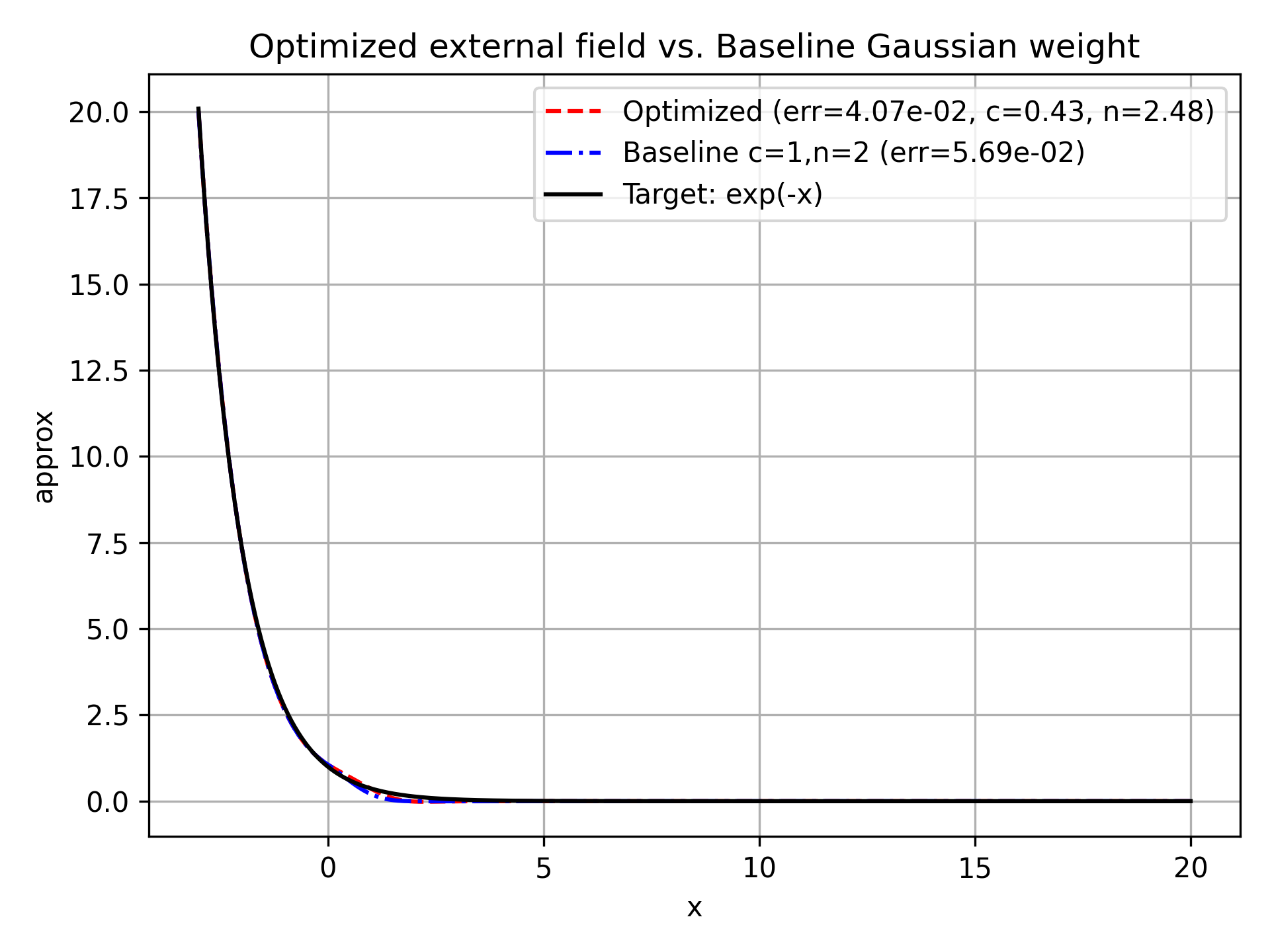}
    \caption{Comparison of optimized and baseline weighted approximants for $e^{-x}$, with learned external field $\Phi(x)=c|x|^{\beta}$ on $x\ge0$.}
    \label{fig:opt_vs_baseline}
\end{figure}

We extend the classical Freud weight $w_\lambda(x)=\exp(-|x|^\lambda)$ by optimizing over the two-parameter family of external fields
\[
  \Phi(x) = c\,|x|^{\beta},\qquad c>0,\;\beta>1,
\]
so that $w(x)=\exp(-\Phi(x))$. (We write $\Phi$ for the external field to avoid any clash with the weighted deep polynomial $Q$, and $\beta$ for the exponent to avoid clash with the polynomial degree $n$.) One checks that $\Phi$ satisfies the hypotheses of Section~\ref{sec:restricted}:
\begin{enumerate}
  \renewcommand{\labelenumi}{(\alph{enumi})}
  \item $\Phi'(x)=c\beta\,x^{\beta-1}>0$ for $x>0$;
  \item $x\Phi'(x)=c\beta\,x^{\beta}$ is strictly increasing on $(0,\infty)$ with $\lim_{x\to0^+}x\Phi'(x)=0$;
  \item \(\displaystyle\frac{x\Phi'(x)}{\Phi(x)}=\frac{c\beta\,x^{\beta}}{c\,x^{\beta}}=\beta\), hence asymptotically constant.
\end{enumerate}
We then fit \((c,\beta)\) \emph{jointly} with the deep polynomial coefficients by minimizing the discrete \(L_2\)-error of the resulting weighted approximant to \(e^{-x}\) on \([-3,20]\). Figure~\ref{fig:opt_vs_baseline} compares the optimized field to the baseline Gaussian choice \((c,\beta)=(1,2)\). Even the simpler choice \((c,\beta)=(1,1)\) yields a nontrivial problem: since the weight acts only on \(x>0\), the growth for \(x<0\) must still be matched by the composite polynomial, whose effective degree grows multiplicatively in the number of layers.

We stress, in response to the concern that the method requires the weight and target to be known analytically, that the weight here is \emph{not} hand-specified: the parameters \((c,\beta)\) are estimated from samples of the target, so the procedure applies whenever one can evaluate the function, even if its growth/decay is not known in closed form. The Airy example of the previous subsection is exactly such a case. A systematic study of joint weight-and-coefficient learning for a broader library of harder targets (asymmetric special functions, profiles with mixed growth/oscillation) is carried out in Section~\ref{sec:finetune}, where a suite of Black--Scholes option-pricing functions with several kinks and one-sided decay is treated at matched coefficient budgets against Chebyshev, Remez, and plain weighted-polynomial baselines.

\section{Fine-tuning weighted deep polynomials with a fixed inner composition}\label{sec:finetune}

\subsection{Motivation}

Recall from Section~\ref{sec:intro} that a deep polynomial with layer degrees $(d_1,\dots,d_L)$ realizes an effective (composite) degree $D=\prod_\ell d_\ell$ from $n_{\mathrm{deep}}=\sum_\ell d_\ell-L+2$ free coefficients. It is this multiplicative growth of $D$ with $L$ that allows the family to outperform a classical polynomial of the same coefficient budget: an informative comparison requires genuine depth, i.e.\ a large $D$.

Training all layers end-to-end to reach a large $D$ is numerically delicate. Under repeated composition, the sensitivity of $P$ to an inner coefficient scales like the product of the derivatives of the outer layers, which grows quickly with $L$; gradient-based training of the full composition correspondingly stalls or diverges, which is why the earlier examples in this paper use only $L=2$. Keeping $L$ small avoids this instability but caps $D$, so the resulting composite is barely more expressive than a single polynomial of equal budget.

We resolve this tension as follows: the inner layers are fixed in advance, chosen from a small dictionary of well-conditioned monotone maps rather than optimized, and only the outer polynomial and the weight are trained. This is analogous to fine-tuning a pretrained network, in which a fixed feature extractor is combined with a trained task-specific head; here the fixed inner composition supplies the effective degree, while the trained outer layer and weight are fit by a convex problem.

\subsection{Construction}\label{subsec:finetune-construction}

\paragraph{Domain conventions.} Given $x\in[a,b]$, let $t=2(x-a)/(b-a)-1\in[-1,1]$ be the standard variable in which Chebyshev polynomials are evaluated, and let $u=(t+1)/2\in[0,1]$. The inner composition acts on $u\in[0,1]$ and returns a value $u'\in[0,1]$, which is mapped back to $z=2u'-1\in[-1,1]$ before the outer Chebyshev polynomial is evaluated. The inner maps below are thus self-maps of $[0,1]$, while the outer polynomial is evaluated on $[-1,1]$ as usual; since each inner map fixes the endpoints, no further rescaling is needed.

\paragraph{Dictionary of inner maps.} We use a small library $\mathcal G$ of monotone increasing polynomials $g:[0,1]\to[0,1]$ with $g(0)=0$, $g(1)=1$, constructed as follows.

For $n\ge1$, let $\phi_n(u)=(2n+1)\binom{2n}{n}u^n(1-u)^n$ on $[0,1]$; this is, up to normalization, the density of a symmetric $\mathrm{Beta}(n+1,n+1)$ random variable, and $\int_0^1\phi_n=1$. Define
\[
p_n(u)\;=\;\int_0^u \phi_n(s)\,ds ,
\]
a polynomial of degree $2n+1$. Since $\phi_n$ vanishes to order $n$ at each endpoint, $p_n(0)=0$, $p_n(1)=1$, and $p_n^{(k)}(0)=p_n^{(k)}(1)=0$ for $1\le k\le n$; by symmetry of $\phi_n$ about $u=1/2$, $p_n(u)+p_n(1-u)=1$. Direct integration gives
\[
p_1(u)=3u^2-2u^3,\qquad p_2(u)=10u^3-15u^4+6u^5,\qquad p_3(u)=35u^4-84u^5+70u^6-20u^7,
\]
the classical cubic, quintic, and septic smoothstep polynomials, each matching one additional derivative at both endpoints as $n$ increases. We include $p_1,p_2,p_3$ in $\mathcal G$.

Because $\phi_n$ is symmetric, each $p_n$ treats both endpoints identically. To obtain maps that concentrate resolution near one endpoint only, we use two further families. First, the unique quadratics with a single vanishing endpoint derivative,
\[
q_L(u)=u^2,\qquad q_R(u)=1-(1-u)^2=2u-u^2,
\]
which are determined by $g(0)=0$, $g(1)=1$, and $g'=0$ at one endpoint only. Second, we compose $p_1$ with itself asymmetrically,
\[
r_L(u)=p_1(u)^2,\qquad r_R(u)=1-\bigl(1-p_1(u)\bigr)^2,
\]
which remain degree-6 polynomials with $g(0)=0$, $g(1)=1$, $g'(0)=g'(1)=0$, but $r_L$ has a higher-order zero at $u=0$ than at $u=1$ (order $4$ versus order $2$), and $r_R$ the reverse. The dictionary is $\mathcal G=\{q_L,q_R,p_1,p_2,p_3,r_L,r_R\}$. Composing several of these, even though each individual map is either symmetric or one-sided, produces inner compositions $\varphi=g_{i_m}\circ\cdots\circ g_{i_1}$ that concentrate resolution asymmetrically, matching the asymmetric transitions in the targets considered below.

\paragraph{Outer fit.}
With the inner composition $\varphi$ fixed, and writing $z_i=2\varphi(u_i)-1$ for the sample points, the outer coefficients $\hat c=(\hat c_k)_{k=0}^{d_{\mathrm{out}}}$ of
\[
Q(x)=w(x;c,\beta,s)\sum_{k=0}^{d_{\mathrm{out}}}\hat c_k\,T_k(z),\qquad
w(x;c,\beta,s)=\exp\Bigl(-c\bigl[\log(1+e^{x-s})\bigr]^\beta\Bigr),
\]
are obtained by solving the discrete minimax problem
\[
\min_{\hat c,\,E}\;E
\quad\text{s.t.}\quad
\Bigl|\,w(x_i;c,\beta,s)\sum_{k=0}^{d_{\mathrm{out}}}\hat c_k T_k(z_i)-f(x_i)\Bigr|\le E,
\qquad i=1,\dots,N .
\]
This is a linear program in $(\hat c,E)$, which we solve using \texttt{scipy.optimize.linprog} with the HiGHS backend. The weight $w$ is the smooth one-sided field of Section~\ref{sec:generalizedweight}, with parameters $c>0$, $\beta>1$, and threshold $s$, and satisfies the hypothesis of Proposition~\ref{prop:compact}. Because the inner composition $\varphi$ enters only through the fixed features $z_i$, the outer fit avoids the ill-conditioning associated with training all layers of a deep composition end-to-end. The weight parameters $(c,\beta,s)$ are refined by a derivative-free search, re-solving the linear program at each trial point. As in Section~\ref{sec:computational}, we stop when the relative change in the objective is below $10^{-12}$, and in the greedy composition step below we accept a new map only when it gives a strict decrease in the discrete uniform error.

\paragraph{Selecting the inner composition.}
We build $\varphi$ greedily. Starting from the identity, at each step we append every $g\in\mathcal G$, solve the outer linear program for each candidate $g\circ\varphi$, and keep the map that most reduces $\|Q-f\|_\infty$ on the sample grid. The candidate is accepted only if it improves the current best error; otherwise the search terminates. This produces a fixed inner feature map and a trained outer weighted Chebyshev expansion.

\paragraph{Cost and comparison protocol.}
Following the cost model of Section~\ref{sec:computational}, we distinguish trainable degrees of freedom from evaluation cost. For the fine-tuned weighted \emph{deep} (composite) polynomial, the trainable parameters are the $d_{\mathrm{out}}+1$ outer Chebyshev coefficients together with the three weight parameters $(c,\beta,s)$, so
\[
\mathrm{dof}=d_{\mathrm{out}}+4 .
\]
The inner composition is fixed after dictionary selection and therefore contributes evaluation cost but no trainable parameters. All methods below are compared at equal trainable degrees of freedom. Thus, at a given budget $\mathrm{dof}$, the Chebyshev and Remez baselines use degree $\mathrm{dof}-1$, while the fine-tuned method uses outer degree $d_{\mathrm{out}}=\mathrm{dof}-4$ plus the learned weight.

The evaluation costs remain comparable, though not identical. At $\mathrm{dof}=64$, the Chebyshev and Remez baselines have degree $63$ and evaluate in about $63$ multiplications by Clenshaw's recurrence. The fine-tuned approximant uses about $60$ multiplications for the outer polynomial, at most roughly $40$ for the inner composition (at most six dictionary maps, each of degree at most $7$, evaluated by Horner's rule), and a constant number of operations for the weight, which is also required by the plain weighted baseline. Hence the composite has a modest constant-factor overhead, less than a factor of two in this setting, while realizing an effective degree of up to
\[
d_{\mathrm{out}}\prod_j \deg g_{i_j},
\]
far beyond the baseline degree. Algorithm~\ref{alg:finetune} summarizes the procedure.

\begin{algorithm}[H]
\caption{Fine-tuned weighted deep polynomial}
\label{alg:finetune}
\DontPrintSemicolon
\KwIn{target $f$, samples $\{x_i\}_{i=1}^N\subset[a,b]$, dictionary $\mathcal G$, outer degree $d_{\mathrm{out}}$, maximum depth $M$, initial weight parameters $(c,\beta,s)$}
\KwOut{inner composition $\varphi$, outer coefficients $\hat c$, weight parameters $(c,\beta,s)$}
$u_i\gets \tfrac12\bigl(2(x_i-a)/(b-a)-1+1\bigr)$\;
\textbf{FitOuter}$(\varphi,c,\beta,s)$: set $z_i\gets 2\varphi(u_i)-1$; solve the LP $\min_{\hat c,E}E$ subject to the minimax constraints above; \Return $(\hat c,E)$\;
$\varphi\gets\mathrm{id}$;\quad $(\hat c,E^\star)\gets \textbf{FitOuter}(\varphi,c,\beta,s)$\;
\For{$m=1$ \KwTo $M$}{
  \ForEach{$g\in\mathcal G$}{ $(\hat c_g,E_g)\gets \textbf{FitOuter}(g\circ\varphi,c,\beta,s)$\; }
  $g^\star\gets\arg\min_{g\in\mathcal G}E_g$\;
  \eIf{$E_{g^\star}<E^\star$}{ $\varphi\gets g^\star\circ\varphi$;\ $E^\star\gets E_{g^\star}$;\ $\hat c\gets \hat c_{g^\star}$\;}{\textbf{break}\;}
}
Refine $(c,\beta,s)$ by derivative-free search, re-running \textbf{FitOuter}$(\varphi,\cdot)$ at each step\;
\KwRet{$\varphi,\ \hat c,\ (c,\beta,s)$}
\end{algorithm}

\subsection{Experiments: Black--Scholes option-pricing functions}\label{subsec:finetune-experiments}

We test the construction on targets generated from the Black--Scholes model \cite{BlackScholes}. If $V=V(S,t)$ denotes the price of a derivative with underlying spot price $S$, time $t$, volatility $\sigma$, and risk-free rate $r$, then $V$ satisfies the Black--Scholes partial differential equation
\[
\frac{\partial V}{\partial t}
+\frac{1}{2}\sigma^2S^2\frac{\partial^2 V}{\partial S^2}
+rS\frac{\partial V}{\partial S}
-rV=0,
\]
with terminal condition determined by the payoff. For a European call with strike $K$ and maturity $T$, the payoff is
\[
V(S,T)=(S-K)_+,
\]
and the corresponding call price is
\[
C(S,K,T,r,\sigma)
=
S\Phi(d_1)-K e^{-rT}\Phi(d_2),
\]
where
\[
d_1=\frac{\log(S/K)+(r+\sigma^2/2)T}{\sigma\sqrt{T}},
\qquad
 d_2=d_1-\sigma\sqrt{T},
\]
and $\Phi$ is the standard normal distribution function. In our experiments the approximation variable is negative log-moneyness,
\[
x=-\log(S/K),
\qquad S=K e^{-x},
\]
so the explicit target used for a vanilla call is
\[
f_{K,T,r,\sigma}(x)
:=
C(K e^{-x},K,T,r,\sigma)
=
K\Bigl(e^{-x}\Phi(d_1(x))-e^{-rT}\Phi(d_2(x))\Bigr),
\]
with
\[
d_1(x)=\frac{-x+(r+\sigma^2/2)T}{\sigma\sqrt{T}},
\qquad
 d_2(x)=d_1(x)-\sigma\sqrt{T}.
\]
In this variable, a call grows like $K e^{|x|}$ as $x\to-\infty$ and decays to zero as $x\to+\infty$, giving exactly the asymmetric growth--decay structure targeted in this paper. At short maturity the payoff also creates a sharp near-kink around the at-the-money region.

The four test functions are built from this explicit call target: a short-maturity call, a call spread, a butterfly spread, and a composite option book. For strikes $K_1<K_2<K_3$, the spread and butterfly targets are
\[
f_{\mathrm{spread}}(x)
=f_{K_1,T,r,\sigma}(x)-f_{K_2,T,r,\sigma}(x),
\]
\[
f_{\mathrm{butterfly}}(x)
=f_{K_1,T,r,\sigma}(x)-2f_{K_2,T,r,\sigma}(x)+f_{K_3,T,r,\sigma}(x),
\]
and the option book is a fixed linear combination
\[
f_{\mathrm{book}}(x)=\sum_{j=1}^m a_j f_{K_j,T_j,r_j,\sigma_j}(x).
\]
Spreads and butterflies superpose several transition regions, while the option book combines multiple maturities and strikes, producing a more challenging multi-scale target.

All targets are sampled on a uniform grid of $N=1200$ points on $[-2,6]$, and all errors are reported on this grid. At each budget $\mathrm{dof}\in\{40,48,64\}$ we compare the fine-tuned weighted \emph{deep} (composite) polynomial against three matched-budget baselines: the Chebyshev polynomial of degree $\mathrm{dof}-1$; the Remez minimax polynomial of the same degree, computed by an exchange iteration in the Chebyshev basis; and a plain weighted polynomial with the same learned weight but no inner composition. The last baseline isolates the contribution of the fixed inner composition.

\begin{figure}[h]
    \centering
    \includegraphics[width=0.8\textwidth]{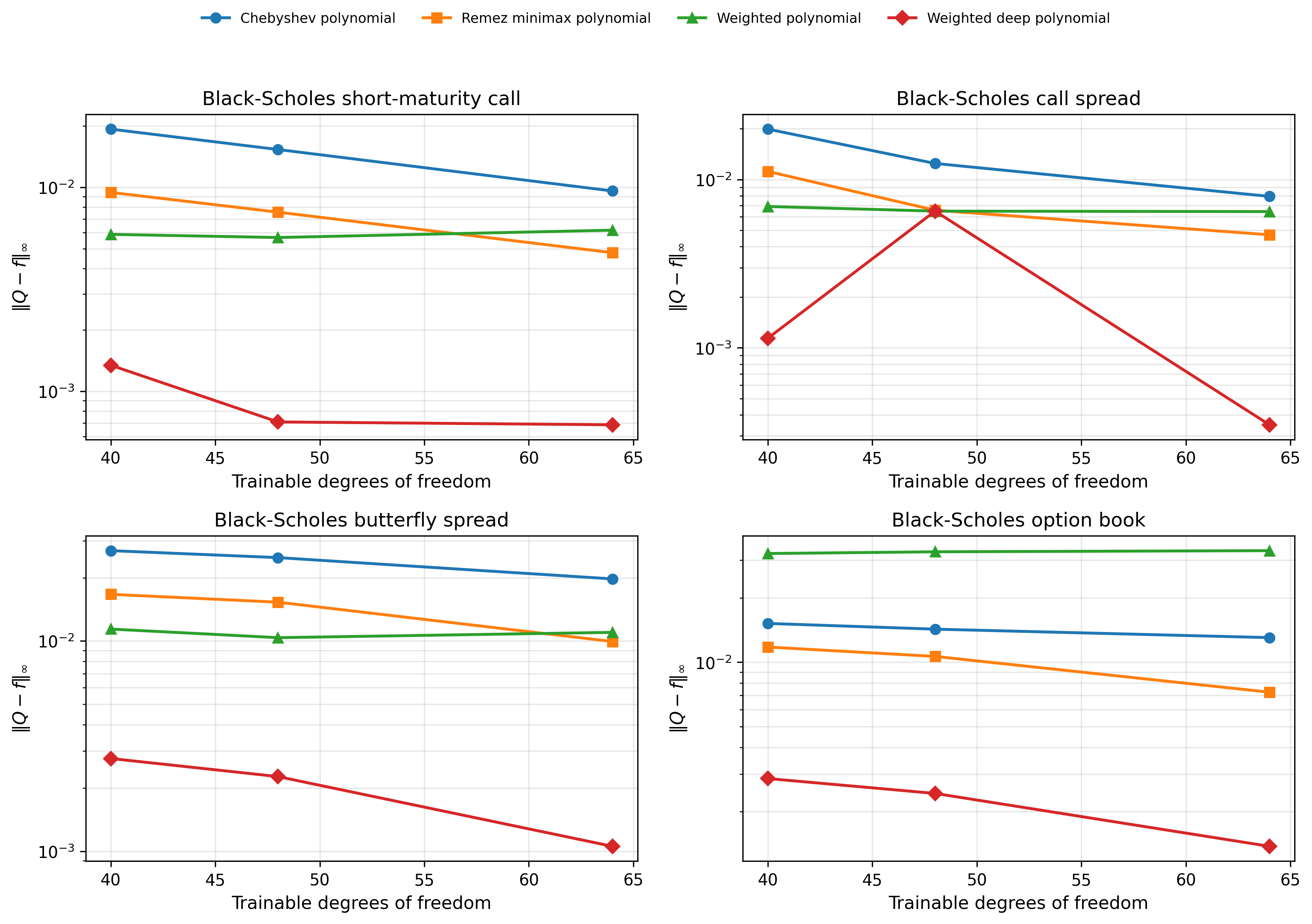}
    \caption{Uniform error $\|Q-f\|_\infty$ versus trainable degrees of freedom for the four Black--Scholes targets, comparing the Chebyshev polynomial, the Remez polynomial, the plain weighted polynomial, and the fine-tuned weighted deep polynomial, all at matched budget.}
    \label{fig:bs_summary_max}
\end{figure}

\begin{figure}[h]
    \centering
    \includegraphics[width=0.8\textwidth]{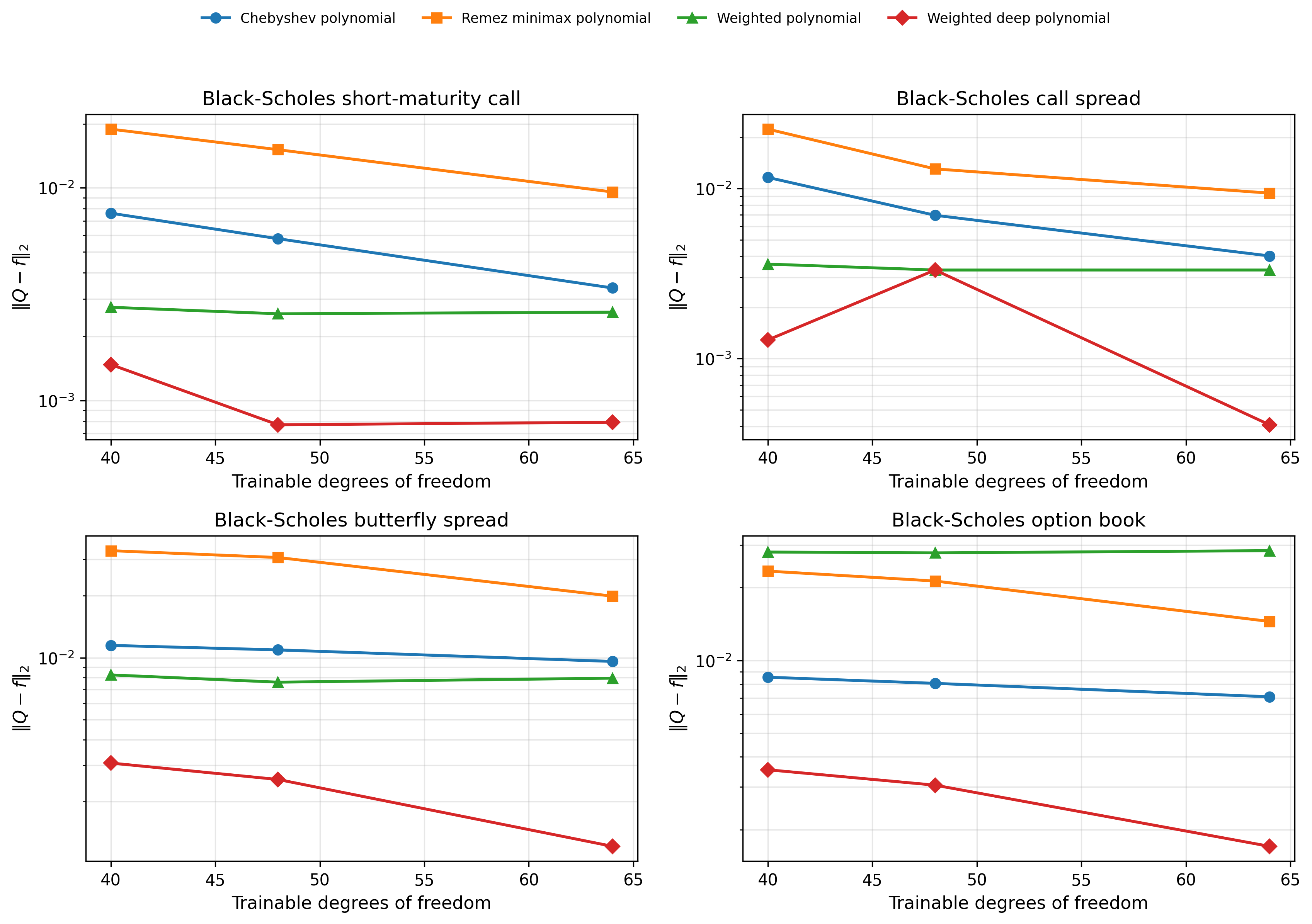}
    \caption{Discrete $L_2$ error $\|Q-f\|_2$ versus trainable degrees of freedom for the same four targets and methods as in Figure~\ref{fig:bs_summary_max}.}
    \label{fig:bs_summary_l2}
\end{figure}

Figures~\ref{fig:bs_summary_max} and \ref{fig:bs_summary_l2} report the uniform and discrete $L_2$ errors as functions of $\mathrm{dof}$. Across all four targets and all three budgets the fine-tuned weighted \emph{deep} (composite) polynomial attains the smallest error, typically by an order of magnitude relative to the Chebyshev and Remez baselines. The plain weighted polynomial is not uniformly competitive: on the option book it performs worse than Chebyshev, so the improvement is attributable to the fixed inner composition rather than the weight alone.

\begin{figure}[h]
    \centering
    \includegraphics[width=0.8\textwidth]{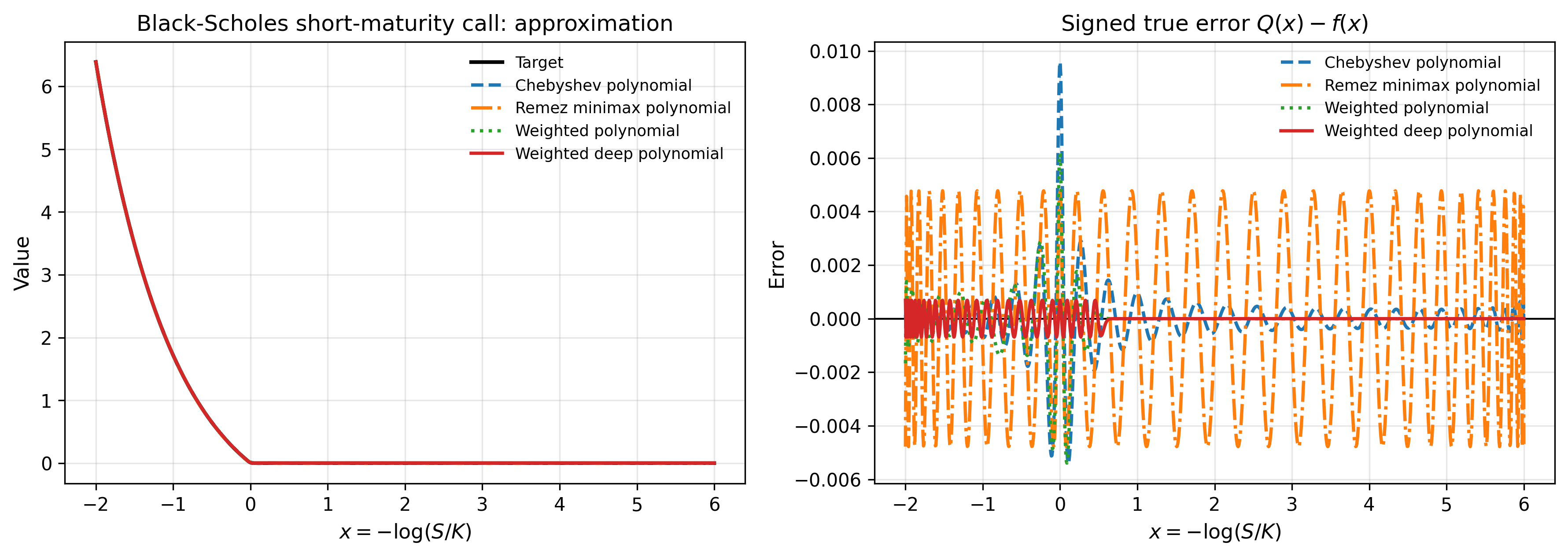}\\[1mm]
    \includegraphics[width=0.8\textwidth]{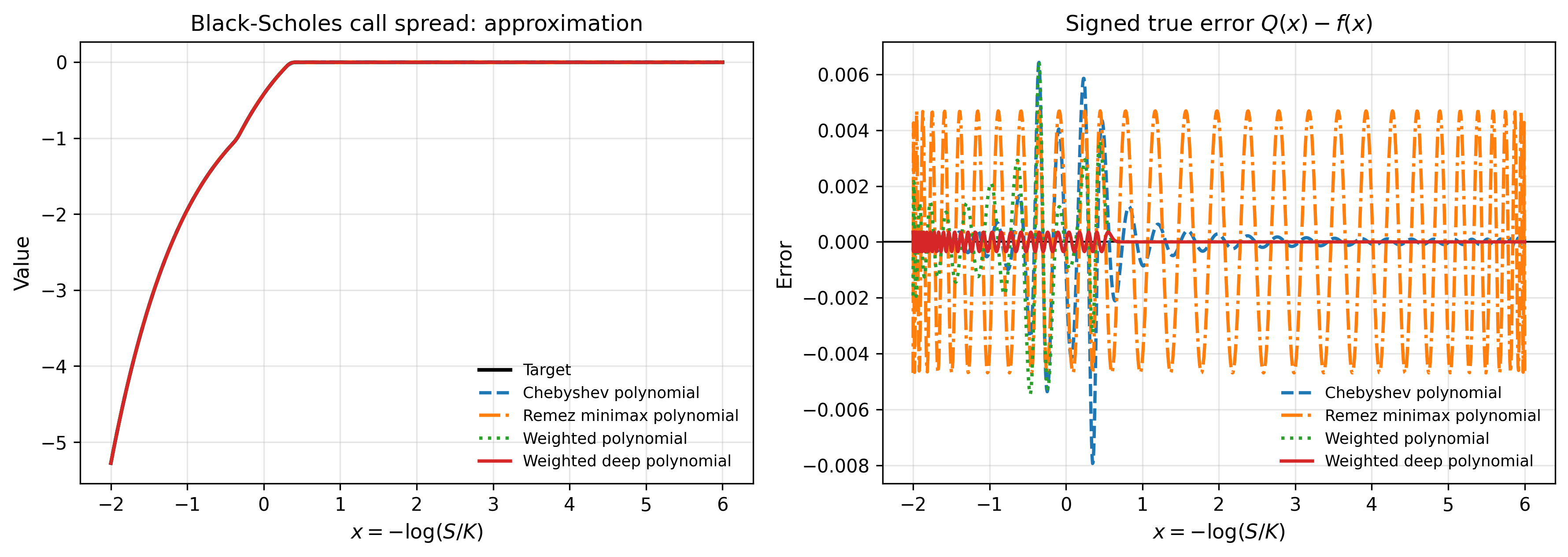}
    \caption{Approximation (left) and signed error $Q(x)-f(x)$ (right) at $\mathrm{dof}=64$ for the short-maturity call (top) and the call spread (bottom).}
    \label{fig:bs_detail_1}
\end{figure}

\begin{figure}[h]
    \centering
    \includegraphics[width=0.8\textwidth]{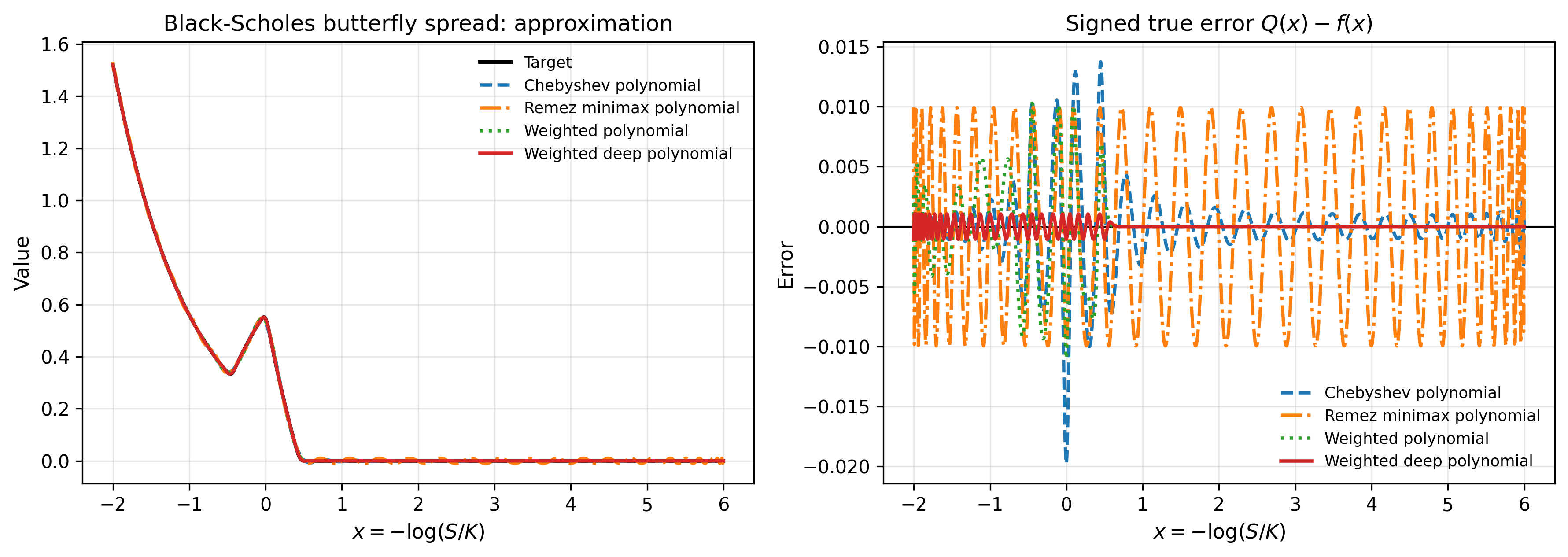}\\[1mm]
    \includegraphics[width=0.8\textwidth]{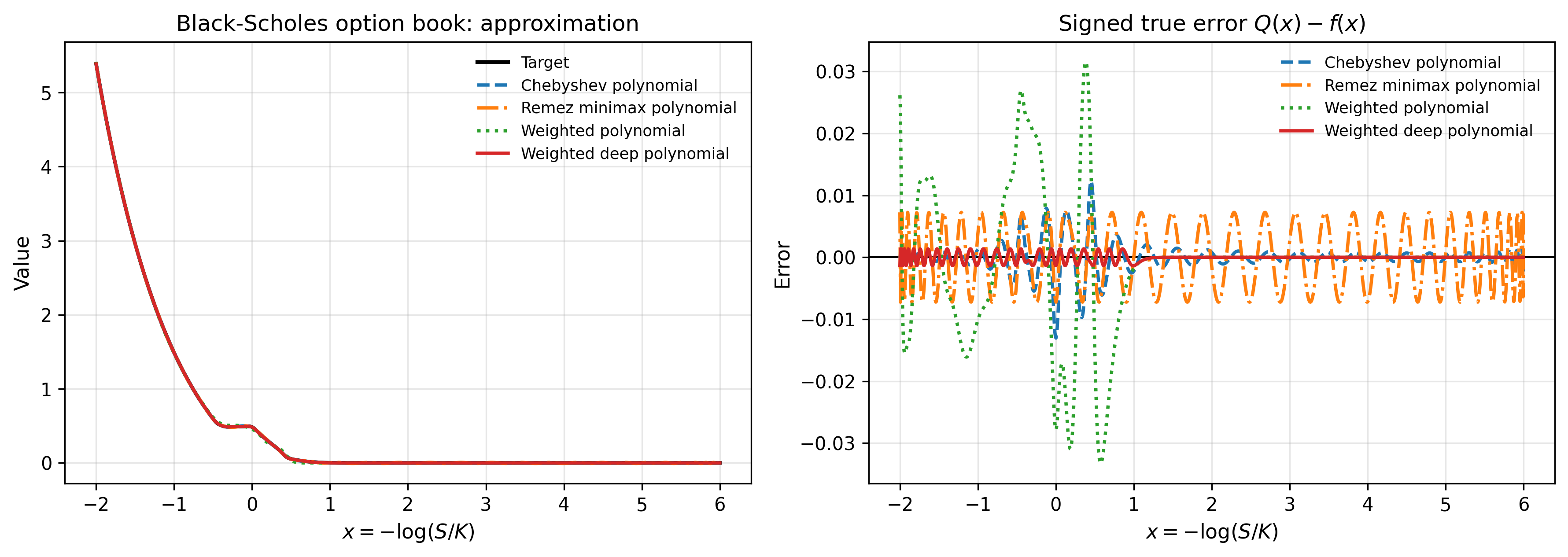}
    \caption{Approximation (left) and signed error (right) at $\mathrm{dof}=64$ for the butterfly spread (top) and the option book (bottom).}
    \label{fig:bs_detail_2}
\end{figure}

Figures~\ref{fig:bs_detail_1} and \ref{fig:bs_detail_2} show the approximation and signed error at $\mathrm{dof}=64$ for all four targets. The Remez polynomial equioscillates, as a minimax fit must, spreading a ripple of order $10^{-3}$ across the entire interval, including the region where the target is exponentially small. The fine-tuned weighted deep polynomial instead concentrates a much smaller error near the at-the-money region and is essentially exact elsewhere. On the option book the plain weighted polynomial develops visible oscillations near the transitions, consistent with its weaker summary errors in Figure~\ref{fig:bs_summary_max}.

For the butterfly spread, we also test the real cosine-sum ESPRIT algorithm of Algorithm~1 in \cite{DPR-ESPRIT}. The method builds a Toeplitz--Hankel matrix from midpoint samples, estimates the frequencies through an SVD-based matrix pencil, and then computes the amplitudes by least squares. In this experiment, ESPRIT gives the largest error among the matched-budget models. This is consistent with the structure of the Black--Scholes butterfly: the function is localized, nonperiodic, nearly flat on much of the interval, and has sharp transition regions near the strikes. A short global cosine sum has to represent both the flat tail and the localized transitions with the same set of frequencies, which makes the pencil estimation sensitive and causes the least-squares amplitude fit to distribute error across the interval. The weighted deep polynomial behaves better here because the one-sided weight captures the decaying side, while the stable inner map reallocates resolution toward the transition region.

\begin{figure}[H]
    \centering
    \begin{minipage}{0.49\textwidth}
        \centering
        \includegraphics[width=\linewidth]{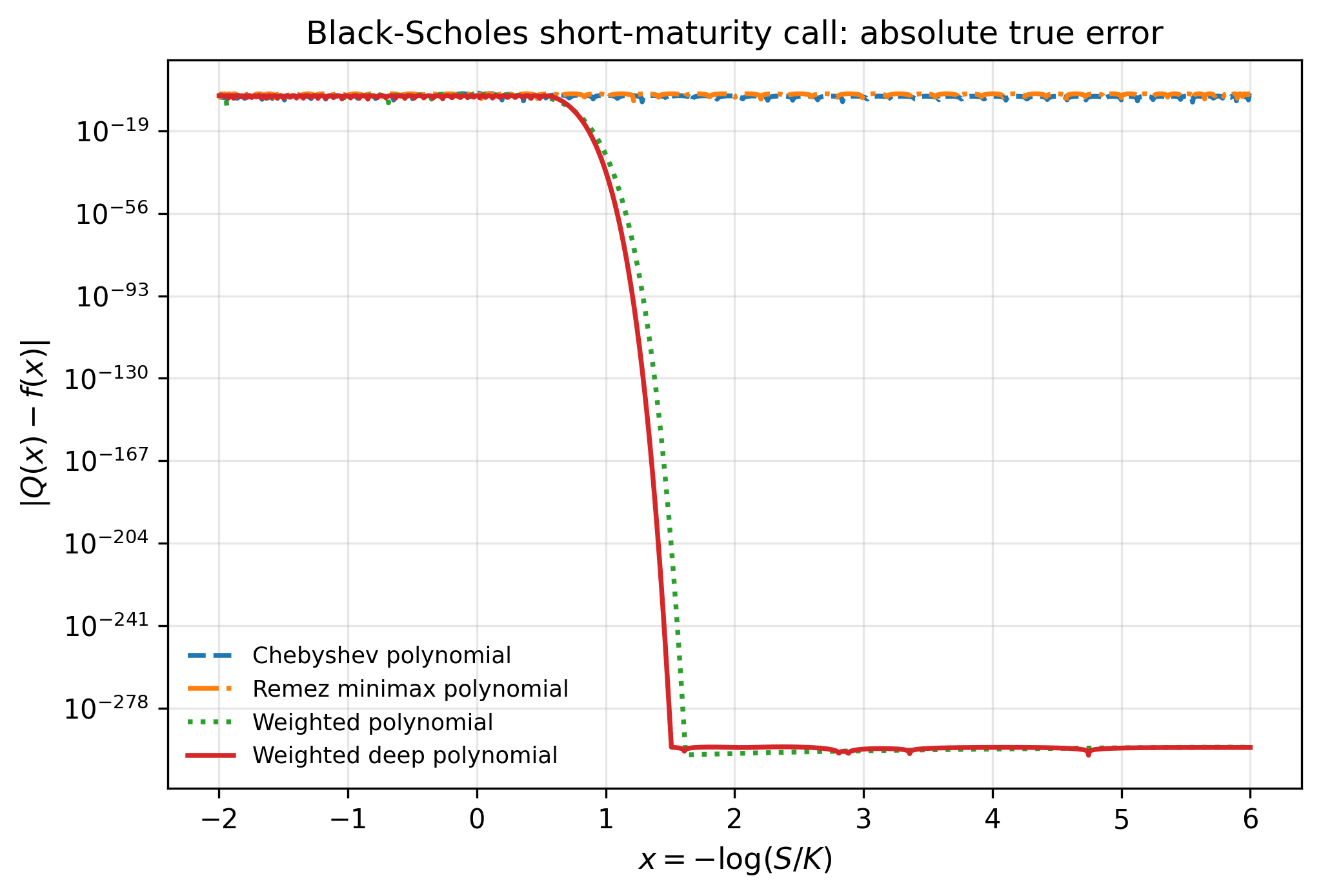}
        \caption*{(a) short-maturity call}
    \end{minipage}\hfill
    \begin{minipage}{0.49\textwidth}
        \centering
        \includegraphics[width=\linewidth]{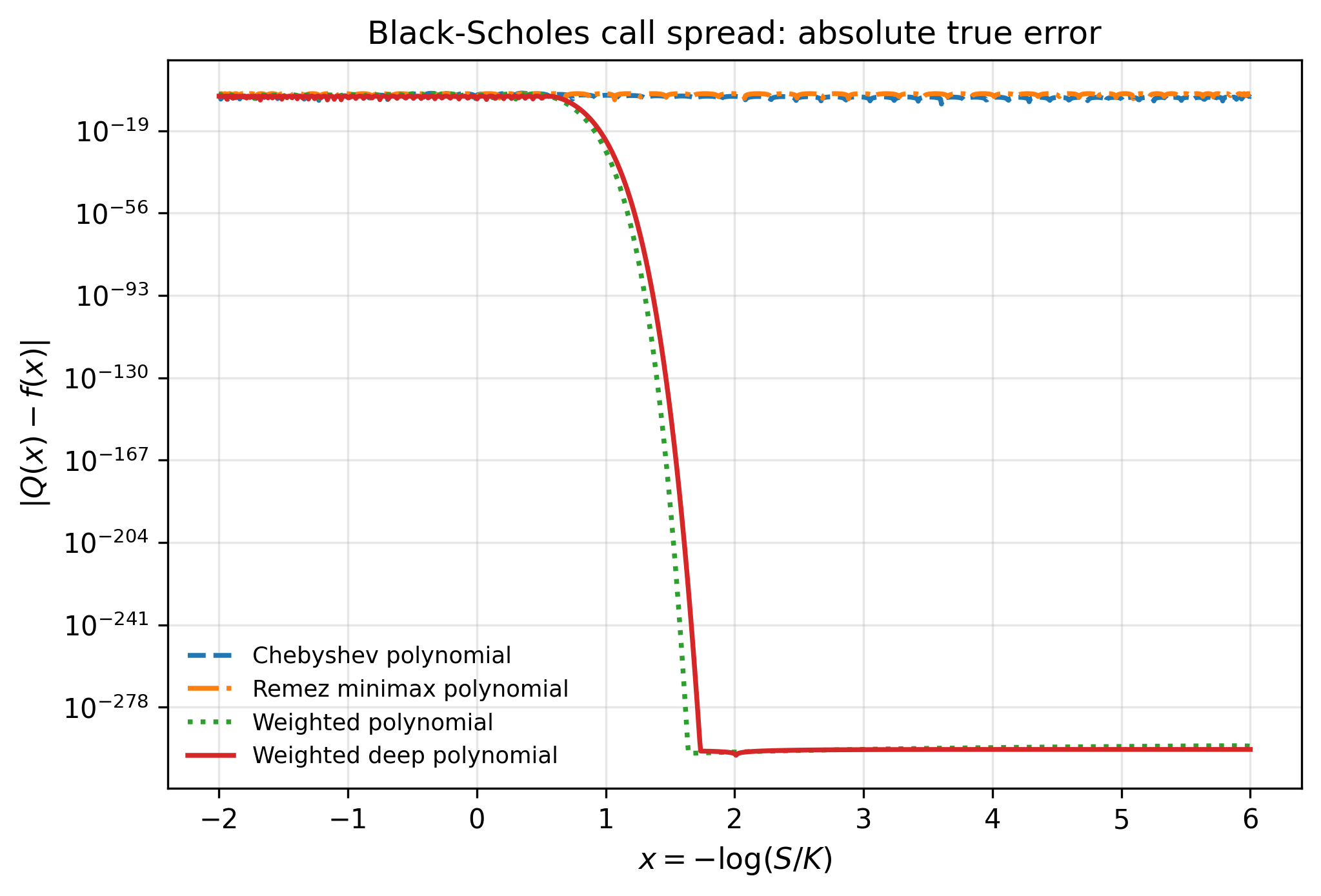}
        \caption*{(b) call spread}
    \end{minipage}

    \vspace{2mm}
    \begin{minipage}{0.49\textwidth}
        \centering
        \includegraphics[width=\linewidth]{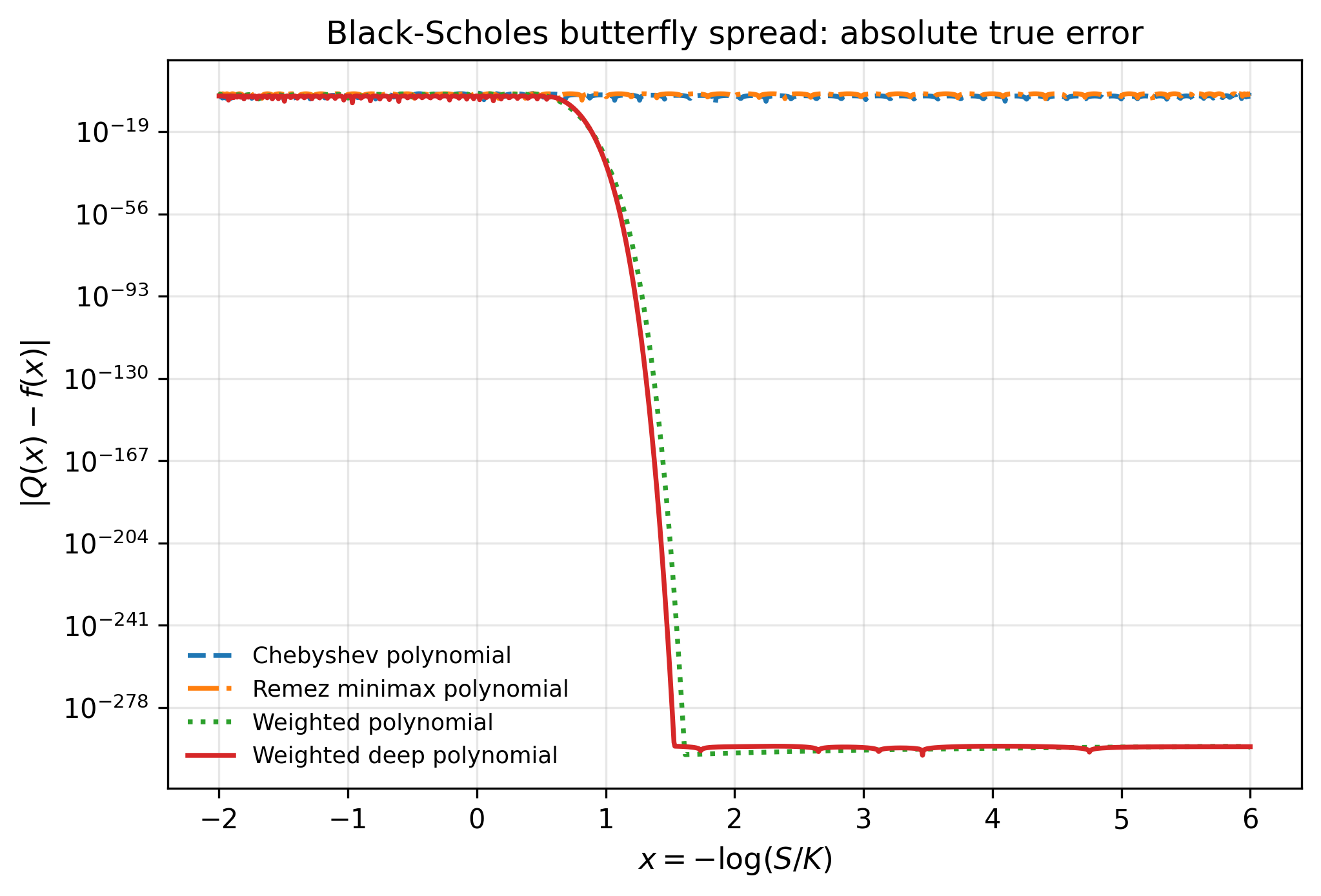}
        \caption*{(c) butterfly spread}
    \end{minipage}\hfill
    \begin{minipage}{0.49\textwidth}
        \centering
        \includegraphics[width=\linewidth]{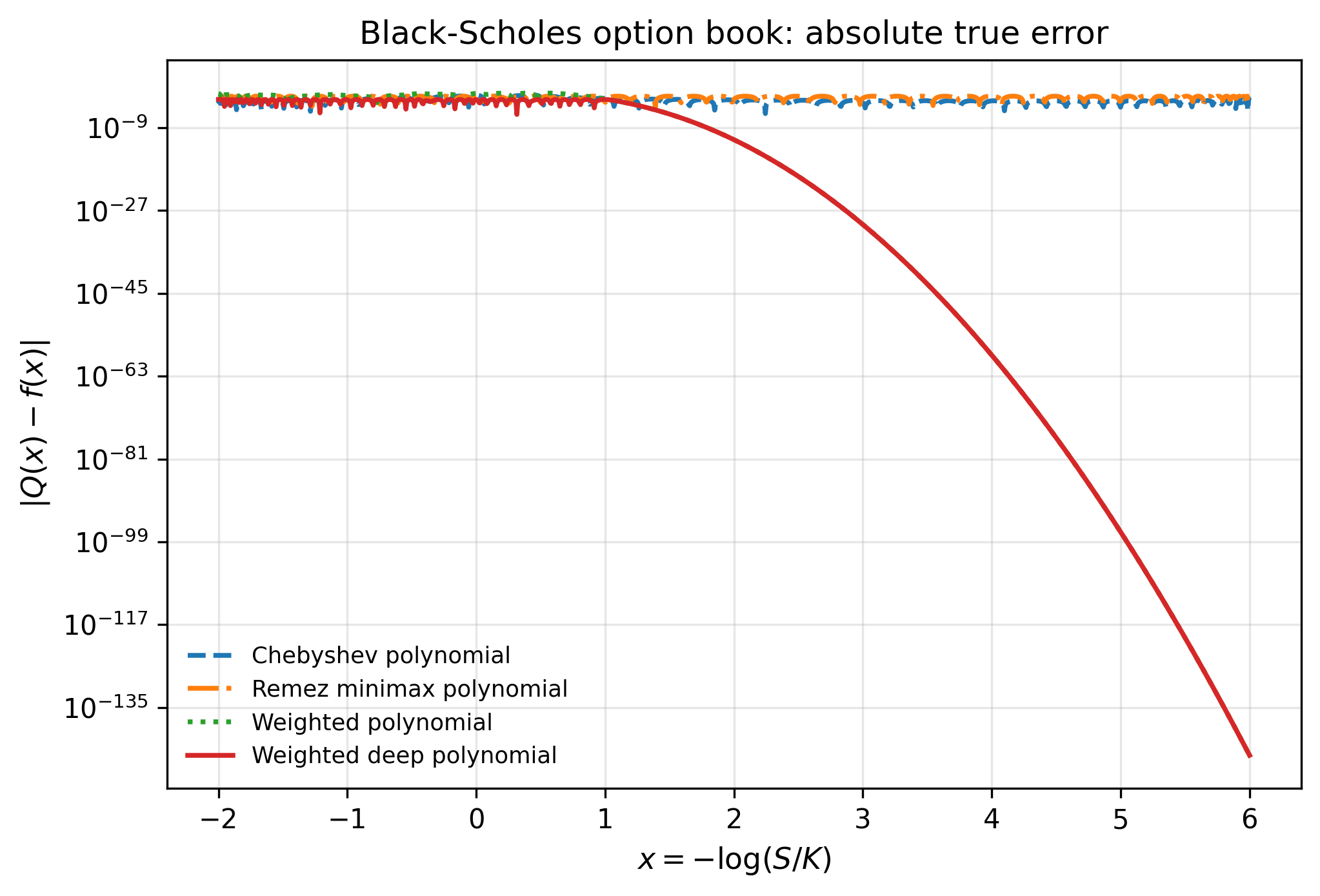}
        \caption*{(d) option book}
    \end{minipage}
    \caption{Absolute error $|Q(x)-f(x)|$ on a logarithmic axis at $\mathrm{dof}=64$. The weighted methods underflow to machine zero on the decaying side, while the Chebyshev and Remez errors remain flat at $O(10^{-3})$. On the option book (d) only the fine-tuned weighted deep polynomial reaches this underflow.}
    \label{fig:bs_abs}
\end{figure}

Figure~\ref{fig:bs_abs} shows absolute error on a logarithmic axis. On the decaying side the one-sided weight suppresses the residual super-exponentially, so the weighted methods underflow to machine zero once $x$ passes the transition, while the Chebyshev and Remez errors remain flat at $O(10^{-3})$ across the domain. The classical polynomial baselines are tail-limited: equioscillation forces them to spend accuracy on the exponentially small tail, leaving less for the kink, whereas the weighted methods pay essentially nothing on the tail. The plain weighted polynomial reaches the same underflow on three of the four targets but not on the option book, where only the fine-tuned construction resolves the tail, again isolating the contribution of the fixed inner composition. This behavior is the numerical counterpart of Proposition~\ref{prop:compact}: beyond the effective interval the weighted approximant contributes an error dominated by the tail of $f$ itself, so the right endpoint of the working interval is immaterial for the weighted methods, and extending $[-2,6]$ further to the right would leave their reported errors unchanged.

These experiments extend the earlier $e^{-x}$ and Airy examples to harder, multi-scale targets and larger budgets, and show that the decaying tail can be resolved to machine precision once the inner composition is fixed and the outer fit is convex. The comparisons reported are against the matched-budget Chebyshev, Remez, and plain weighted-polynomial baselines shown. The construction recovers the expressive power of deep composition while avoiding the instability of end-to-end training, at an evaluation cost of the same order as the classical baselines.






\end{document}